\documentclass[10pt]{amsart}
\usepackage{amssymb, mathrsfs,amsmath}
\usepackage{latexsym}
\usepackage{amsfonts}
\usepackage{shadow}

\newtheorem{Pa}{Paper}[section]
\newtheorem{Tm}[Pa]{{\bf Theorem}}
\newtheorem{La}[Pa]{{\bf Lemma}}
\newtheorem{Cy}[Pa]{{\bf Corollary}}
\newtheorem{Ob}[Pa]{{\bf Observation}}

\newtheorem{Pn}[Pa]{{\bf Proposition}}

\newtheorem{Ex}[Pa]{{\bf Example}}

\makeatletter
\def\Ddots{\mathinner{\mkern1mu\raise\p@
\vbox{\kern7\p@\hbox{.}}\mkern2mu
\raise4\p@\hbox{.}\mkern2mu\raise7\p@\hbox{.}\mkern1mu}}
\makeatother
\def\C{\mathbb C}

\date{}
\author[D. Alpay]{Daniel Alpay}
\author[P. Jorgensen]{Palle Jorgensen}
\address{(DA) Department of Mathematics
\newline
Ben Gurion University of the Negev \newline P.O.B. 653,
\newline
Be'er Sheva 84105, \newline ISRAEL} \email{dany@math.bgu.ac.il}
\address{(PJ)
Department of Mathematics\newline 14 MLH \newline The University
of Iowa, Iowa City,\newline IA 52242-1419 USA}
\email{jorgen@math.uiowa.edu}
\author[I. Lewkowicz]{Izchak Lewkowicz}
\address{(IL) Department of Electrical Engineering
\newline
Ben Gurion University of the Negev \newline P.O.B. 653,
\newline
Be'er Sheva 84105, \newline ISRAEL}
\email{izchak@ee.bgu.ac.il}
\thanks{This research is partially supported by the BSF
grant no. 2010117.}
\thanks{D. Alpay thanks the
Earl Katz family for endowing the chair which supported his
research.}
\title
[
families of rectangular, para-unitary FIR systems
]
{
characterizations of families of rectangular,\\[0.2cm]
finite impulse response, para-unitary systems}
\begin{document}
\begin{abstract}
We here study Finite Impulse Response (FIR) rectangular, not
necessarily causal, systems which are (para)-unitary on the unit
circle (=the class $\mathcal{U}$). First, we offer three
characterizations of these systems. Then, introduce a
description of all FIRs in $\mathcal{U}$, as copies of a real
polytope,
parametrized by the dimensions and the McMillan degree
of the FIRs.

Finally, we present six simple ways (along with their combinations)
to construct, from any FIR, a large family of FIRs, of various
dimensions and McMillan degrees, so that whenever the original system
is in $\mathcal{U}$, so is the whole family.

A key role is played by Hankel matrices.
\end{abstract}

\noindent
 \maketitle

\section{Introduction}
\setcounter{equation}{0}
\label{sec:intoduction}

This work is on the crossroads of Operator and Systems theory
from the mathematical side and Control, Signal Processing and
Communications theory from the engineering side. It addresses
problems or employs tools from all these areas. Thus, it is
meant to serve as a bridge between the corresponding communities.
We start by formally laying out the set-up.

\subsection{Finite Impulse Response}

We here focus on $p\times m$-valued polynomials
of a complex variable $z$, of the form
\begin{equation}\label{eq:poly}
F(z)=z^{q}\left(z^{-1}B_1+~\ldots~+z^{-n}B_n\right),
\end{equation}
where the natural $n$ and the integer $~q~$ are parameters.
Hereafter, we relate to Laurent polynomials when the powers
of $z$ may be positive or negative (or both).
\vskip 0.2cm

The polynomial $F(z)$ in \eqref{eq:poly} may be viewed as the
(two-sided) $Z$-transform of $~\Phi(t)$ a (discrete) time sequence,
\begin{equation}\label{eq:time}
\Phi(t)=
\delta_K(t-q+1)B_1+~\ldots~+\delta_K(t-q+n)B_n
\quad\quad\quad t\quad{\rm integral~~ variable},
\end{equation}
where $\delta_K$ is the Kronecker delta,
\[
\delta_K(\beta)=\left\{\begin{smallmatrix}1&~&\beta=0,\\~\\
0&~&\beta\not=0,\end{smallmatrix}\right.
\quad\quad\quad\quad\beta\quad{\rm integer}.
\]
Hence, in Engineering terminology $\Phi(t)$ (and often $F(z)$) is
referred to as Finite Impulse Response (i.e. the support of $\Phi(t)$ is
finite).
\vskip 0.2cm

Moreover $F(z)$ will be called {\em causal} whenever $~1\geq q~$ (i.e.
$\Phi(t)\equiv 0$ for all $0>t$) and {\em strictly  causal} if $~0\geq q$.
Similarly, {\em (strictly) anti-causal} when ($q\geq n+1$) $q\geq n$.
When anti-causal, $F(z)$ in \eqref{eq:poly} is a usual polynomial with
non-negative powers of the variable $~z$.
\vskip 0.2cm

Finite Impulse Response functions (=Laurent polynomials)
are of numerous applications in communications control and
signal processing see e.g. \cite{CB1}, \cite{CB2}, \cite{CZZM3},
\cite{Dav}, \cite{LM1}, \cite{LM2}, \cite{TV}, \cite{VHEK}.

\subsection{Unitary symmetry}
\label{subsec:UnitarySym}

Let ${\mathbb T}$ be the unit circle,
\[
{\mathbb T}:=\{z\in\C~:~|z|=1~\}.
\]
We denote by $\mathcal{U}$ the class of $p\times m$-valued
rational functions having unitary symmetry on the unit circle, i.e.
\begin{equation}\label{eq:DefU}
\mathcal{U}:=\left\{~F(z)~:~\left\{
\begin{matrix}\left(F(z)\right)^*F(z)\equiv
I_m&p\geq m&{\rm isometry}\\~\\
F(z)\left(F(z)\right)^*\equiv I_p&m\geq p&
{\rm co-isometry}\end{matrix}\right.
\quad\quad\forall z\in{\mathbb T}\right\}.
\end{equation}
In engineering terminology, if in addition all poles of $F(z)$ are
within the open unit disk (=Schur stable), the function $F(z)$ is
called ~{\em lossless}\begin{footnote}{Passive electrical
circuits are either dissipative or lossless.}\end{footnote},
see e.g. \cite{GVKDM}, \cite[Section 14.2]{Vaid}
or ~{\em all-pass}\begin{footnote}{For example, in studying
classical filters a ``high-pass" could be viewed as an
``all-pass" minus a ``low-pass".}\end{footnote}.
\vskip 0.2cm

For a given $p\times m$-valued rational function $F(z)$ we here
define the \mbox{$m\times p$-valued} ~{\em conjugate}~ 
rational function as,
\[
F^{\#}(z):=\left(F\left(\frac{1}{z^*}\right)\right)^*.
\]
Note that
\[
{F^{\#}(z)}_{|_{z\in\mathbb{T}}}=
\left({F(z)}_{|_{z\in\mathbb{T}}}\right)^*.
\]
It is well known,  see e.g. \cite[Eq. (3.1)]{AlGo1},
\cite[D\'{e}finition 35]{Icart}, that for rational functions condition
\eqref{eq:DefU} is equivalent to the following,
\[
F\in\mathcal{U}
\quad\quad\Longleftrightarrow\quad\quad\quad\forall z\in\C\quad
\left\{
\begin{matrix}F^{\#}(z)F(z)\equiv I_m&p\geq m&{\rm isometry},\\~\\
F(z)F^{\#}(z)\equiv I_p&m\geq p&{\rm co-isometry}.\end{matrix}\right.
\]
The interest in the class $\mathcal{U}$ stems from 
various aspects: For realization see
e.g. \cite{AlGo1}, \cite{AlGo2}, 
\cite{AlRak2}
\cite{GVKDM} 
and \cite{Valli} for factorization see e.g.
\cite{Pot} and for some signal
processing applications see e.g. \cite{BJ1}.
\vskip 0.2cm

\subsection{The current work}
The interest in Finite Impulse Response functions within
${\mathcal U}$ (=para-unitary, in signal processing
``dialect") is vast, see e.g.  the books \cite{BJ2},
\cite{J3}, \cite[Section 7.3]{Ma}, \cite[Section 5.2]{SN},
\cite[Section 6.5]{Vaid}, the theses
\cite{Icart}, \cite{Oli} and the papers 
\cite{AJL3}, \cite{AJL4}, 
\cite{AJLM}, \cite{BJ1}, \cite{GNS},
\cite{J2}, \cite{J5}, \cite{OTHN}, \cite{RMW}, 
\cite{SDeLID}
and \cite{YaZh08}.
\vskip 0.2cm

For example, the classical spectral factorization of self adjoint
matrices and the singular values decomposition of rectangular
matrices (for constant matrices see e.g.
\cite[Theorems 2.5.6, 4.1.5]{HJ1} and \cite[Theorem 7.3.5]{HJ1},
respectively) have been generalized to
matrix valued polynomials (where at least one of the factors
is para-unitary). These extension have several signal processing
applications and were studied in \cite{MWBCR} and \cite{TAR}.
\vskip 0.2cm

Similarly, the Q-R factorization of constant matrices has been
generalized to matrix-valued polynomials (one of which is
para-unitary). In \cite{CB1}, \cite{CB2} this was studied along
with applications to OFDM communications
\vskip 0.2cm

{\em
This work focuses on characterizations of families of rectangular
(not necessarily causal) Finite Impulse Response (FIR) functions
within ${\mathcal U}$.
}
\vskip 0.2cm

This work is aimed at three different communities: mathematicians
interested in classical analysis, signal processing engineers and
system and control engineers. Thus adopting the terminology
familiar to one audience, may intimidate or even alienate the
other. For example, as already mentioned, matrix-valued Laurent
polynomials (powers of various signs) and not necessarily causal
Finite Impulse Response systems, are virtually the same entity
seen by a different community.
Similarly, what is known to engineers as McMillan degree also
arises in geometry of loop groups as an index.
\vskip 0.2cm

Books like 
\cite{BJ2}, \cite{HSK}, \cite{SN}, and the
theses \cite{Icart}, \cite{Oli} have made an effort to be at least
``bi-lingual". Lack of space prevents us from providing even a
concise dictionary of relevant terms. Thus, we shall try to employ
only basic concepts.
\vskip 0.2cm

This work is organized as follows.
\vskip 0.2cm

In Subsection \ref{sec:realization1} we construct the Hankel
matrix $\mathbf{H}$ associated with a polynomial $F(z)$ in
\eqref{eq:poly} and show how that McMillan degree can be extracted
from it. In Subsection \ref{subsec:families} we show how to
construct from a single Laurent polynomial a whole family of
Laurent polynomials of various powers and degrees. Moreover,
whenever the original polynomial is in $\mathcal{U}$, so is the
whole resulting family. This construction is based on the Hankel
matrix $\mathbf{H}$. The details are relegated to the Appendix.
\vskip 0.2cm

In subsection \ref{sec:RealizeMatrix} we present a characterization
of causal Schur stable rational functions in $\mathcal{U}$ through
their minimal realization matrix.
In subsection \ref{sec:Blaschke} we present  the Blaschke-Potapov
characterization of rational functions in $\mathcal{U}$.
Here, the rational functions may have poles everywhere (including
infinity) except the unit circle. In Theorem \ref{Tm:UnitaryFIR}
we introduce three convenient formulations of this result, in
the framework of FIRs.
\vskip 0.2cm

As a by-product, this enables us to offer an easy-to-use
description of all FIRs in $\mathcal{U}$ with McMillan
degree and dimensions as parameters. In fact, for causal
systems, it turns out to be a genuine real polytope and in
general these are copies of this polytope. This is in particular
convenient if one wishes to:
(i) Design FIRs within $\mathcal{U}$ through optimization,
see e.g. \cite{RMW}.\\
(ii) Iteratively apply para-unitary similarity, see e.g.
\cite[Section 3.3]{Icart}, \cite{MWBCR}, \cite{SDeLID}.
In signal processing literature, this is associated with
with {\em channel equalization}~ and in communications
literature with ~{\em decorrelation of signals}~ or
(iii) Iteratively apply Q-R factorization in the
framework of communications, see e.g. \cite{CB1}, \cite{CB2}.
\vskip 0.2cm
\vskip 0.2cm

In Subsection \ref{subsec:IsometricHankel} we return to the Hankel
$\mathbf{H}$ and show that more information is ``encoded" in it: Out of
$\mathbf{H}$ one can check whether or not $F(z)$ is in
$\mathcal{U}$ (see Theorem \ref{Ob:UnitaryF}). Moreover, just
from the singular
values of $\mathbf{H}$ one can deduce whether $F(z)$ is square or
rectangular (see Proposition \ref{Pn:UnitaryF}). It turns out
that the conditions there are closely related to those in the
Nehari problem where one approximates an anti-causal rational
polynomial, by a causal one, see e.g. \cite[Section 12.8]{HSK}.
\vskip 0.2cm

A concluding remark is given at the end.

\section{the Hankel matrix and FIR}
\setcounter{equation}{0}
\label{sec:HankelFIR}

In this section we present our first use of Hankel matrices.
They will be re-appear in Section \ref{sec:Hankelrevisited}
and the Appendix.

\subsection{Realization of FIR and the Hankel matrix}
\label{sec:realization1}

In this subsection we focus on the McMillan degree a
$p\times m$-valued polynomial polynomial $F(z)$
\begin{equation}\label{eq:Poly}
F(z)=z^{q}\left(z^{-1}B_1+~\ldots~+z^{-n}B_n\right),
\end{equation}
where the natural $n$ and the integer $~q~$ are parameters.
\vskip 0.2cm

We here specialize textbook material on state space
realization (of not necessarily causal systems) and the
corresponding Hankel matrices, see e.g.
\cite[Subsection 12.8.1]{HSK}.
\vskip 0.2cm

We start by recalling some classical facts concerning realization
of $p\times m$-valued rational function $F(z)$.
With a slight abuse of notation,
an arbitrary $p\times m$-valued rational function $F(z)$ may be
written as
\begin{equation}\label{eq:rational}
F(z)=F_l(z)+D+F_r(z)
\end{equation}
where $D$ is a constant matrix and
\[
\lim\limits_{z\rightarrow ~0}F_l(z)=0_{p\times m}
\quad\quad\quad\quad
\lim\limits_{z\rightarrow\infty}F_r(z)=0_{p\times m}~,
\]
(the subscripts stand for ``left" and ``right"). Note that $F_l(z)$,
$F_r(z)$ may be viewed as the (two-sided) $Z$-transform of
{\em strictly
anti-causal},~ {\em strictly causal}~ (discrete) time sequences,
respectively.
\vskip 0.2cm

Recall that $\tilde{F}_r(z)$ a $p\times m$-valued rational function
with no poles at infinity\begin{footnote}{=Bounded at infinity and
in engineering called ~{\em causal}~ or colloquially ~{\em 
proper}.}\end{footnote} is given by,
\begin{equation}\label{eq:RealizeF}
\tilde{F}_r(z):=F_r(z)+D=C(zI_n-A)^{-1}B+D.
\end{equation}
Sometimes it is convenient to present $\tilde{F}_r(z)$ in 
\eqref{eq:RealizeF}
by its $(\nu+p)\times(\nu+m)$ realization matrix $R$, i.e.
\begin{equation}\label{eq:R}
R:={\footnotesize
\left(\begin{array}{c|c}A&B \\ \hline C&D
\end{array}\right)}.
\end{equation}
A realization of $\tilde{F}_r(z)$ in \eqref{eq:RealizeF}, is called
{\em minimal} if in
\eqref{eq:R}, $\nu$ the
dimension of $A$, is the smallest possible. This $\nu$ is called
the ~{\em McMillan degree}~ of $\tilde{F}_r(z)$ in \eqref{eq:RealizeF}.
\vskip 0.2cm

As already mentioned, in the polynomial framework in \eqref{eq:Poly},
$F(z)$ is causal when $1\geq q$ and strictly causal when $0\geq q$. A
special attention will be devoted 
to the case where in \eqref{eq:Poly} $q=0$, i.e.
\begin{equation}\label{eq:HatFr}
F_o(z):={\tilde{F}_r(z)}_{|_{q=0}}=F(z)_{|_{q=0}}=z^{-1}B_1+~\ldots~,~z^{-n}B_n~,
\end{equation}
\vskip 0.2cm

For future reference, we introduce the following notation\begin{footnote}{In
the sequel, boldface
characters will stand for block-structured matrices.}\end{footnote}
\begin{equation}\label{eq:B_eta}
{\mathbf B}_{\mathbf \eta}:=\left(\begin{smallmatrix}
0_{\eta\cdot p\times m}\\ B_1\\ \vdots\\~\\B_n\end{smallmatrix}\right)
\quad\quad\quad\quad\eta=0,~1,~2,~\ldots
\end{equation}
and by ${\mathbf J}_k$
the following $kp\times kp$ block-shift matrix,
\[
{\mathbf J}_k=\left(\begin{smallmatrix} 
~            &I_p& ~    &~  \\
~            &~  &\ddots&~  \\
~            & ~ & ~    &I_p\\
0_{p\times p}& ~ & ~    &~
\end{smallmatrix}\right).
\]
With this notation, specializing the realization matrix \eqref{eq:R}
to the case of $F_o(z)$ in \eqref{eq:HatFr}, yields
\begin{equation}\label{eq:StraightforwardRealiz}
R:={\footnotesize
\left(\begin{array}{c|c}{\mathbf J}_n&{\mathbf B}_0
\\ \hline I_p\quad 0_{p\times(n-1)p}&0_{p\times m}
\end{array}\right)},
\end{equation}
where ${\mathbf B}_0$ is as in \eqref{eq:B_eta} (with $\eta=0$).
\vskip 0.2cm

This classical approach has two limitations:\\
(i) Often the realization in \eqref{eq:StraightforwardRealiz} 
is not minimal. Moreover, the actual McMillan degree (bounded
from above by $np$) is not apparent from \eqref{eq:StraightforwardRealiz}.
\vskip 0.2cm

(ii) Strictly speaking \eqref{eq:R} and its special case in
\eqref{eq:StraightforwardRealiz} are realization around zero,
suitable for causal systems. Realizations for the anti-causal case
(around $z=\infty$), are fairly common as well.
%
Although known, it is more challenging to write down realizations
of polynomials $F(z)$ containing powers of mixed signs, i.e. when
in \eqref{eq:Poly} $q\in[2,~n-1]$.
\vskip 0.2cm

We next show that representing realizations of FIR systems through
Hankel matrices\begin{footnote}{In general, the Hankel operator is
infinite, but
since we here focus on $F(z)$ in \eqref{eq:rational} with a
Finite Impulse Response, the corresponding Hankel matrix is finite and
no truncation is needed.}\end{footnote}, circumvent both limitations:
\vskip 0.2cm

It is suitable for
all $q$ in \eqref{eq:Poly} and the McMillan degree is apparent.
\vskip 0.2cm

In Subsection \ref{subsec:IsometricHankel} below these Hankel matrices
will be used to introduce a characterization of Finite Impulse
Response functions within $~{\mathcal U}$.
\vskip 0.2cm

To study the Hankel matrix representation of $F(z)$ in \eqref{eq:Poly},
\[
F(z)=z^{q}\left(z^{-1}B_1+~\ldots~+z^{-n}B_n\right),
\quad\quad\quad q\quad{\rm integral~parameter},
\]
we find it convenient to separately consider two
extremes possibilities, and an intermediate case:
\vskip 0.2cm

(i) $q\geq n+1~$ so $F(z)$ is strictly anti-causal and
\vskip 0.2cm

(ii) $0\geq q~$ so $F(z)$ is strictly causal,
\vskip 0.2cm

(iii) $q\in[1,~n]~$ so $F(z)$ is a genuine Laurent polynomial.
\vskip 0.2cm

(i) Assume now that $F(z)$ in \eqref{eq:Poly} is 
{\em strictly causal},~
$-q:=\eta\geq 0$. Here \eqref{eq:Poly} takes the form
\[
F(z)=F_r(z)=z^{-(1+\eta)}B_1+~\ldots~+z^{-(n+\eta)}B_n
\quad\quad\quad\quad\eta\geq 0.
\]
We shall denote by ${\mathbf H_{\eta}}$ (recall, block-structured
matrices are represented by boldface characters) the associated
$p(n+\eta)\times m(n+\eta)~$ Hankel matrix,
\begin{equation}\label{eq:ContrR}
{\mathbf H}_{\eta}
=\left(\begin{smallmatrix}
~     &~     &B_1   &\cdots&B_n\\
~     &\Ddots&   ~  &\Ddots&   \\
B_1   &\cdots&B_n   & ~    & ~ \\
\vdots&\Ddots&~     & ~    & ~ \\
B_n&~ & ~    & ~    & ~    & ~~
\end{smallmatrix}\right).
\end{equation}
The $p\times m$ (block) elements of ${\mathbf H_{\eta}}$ are known
as the Markov parameters of $F_r(z)$ and in particular, 
the first (block) row of ${\mathbf H_{\eta}}$ is the impulse
response of $F_r(z)$.
\vskip 0.2cm

For completeness we add that the Hankel matrix can be obtained
as a product of the observability and controllability matrices
here it takes the form
\[
{\mathbf H}_{\eta}=
\left(\begin{matrix}\mathbf{J}^0&~&\mathbf{J}^1&~&\cdots&~&
\mathbf{J}^{n+\eta-1}
\end{matrix}\right)
\left(\begin{smallmatrix}
{\mathbf B}_{\mathbf \eta}&~     &~                         \\
~                         &\ddots&~                         \\
~                         &~     &{\mathbf B}_{\mathbf \eta}
\end{smallmatrix}\right),
\]
where $~{\mathbf B}_{\mathbf \eta}~$ is as in \eqref{eq:B_eta},
$~\mathbf{J}=\mathbf{J}_{n+\eta}$
and $~\mathbf{J}^0=I_{(n+\eta)p}~$.
\vskip 0.2cm

Associating the Hankel matrix ${\mathbf H}_{\eta}$ in \eqref{eq:ContrR}
with the polynomial $F_r(z)$ is fairly classical and goes back
at least to \cite[Eq. (7)]{GLR}.
\vskip 0.2cm

(ii) The other extreme is where $F(z)$ is {\em strictly 
anti-causal}, i.e. in \eqref{eq:Poly} $q\geq n+1$ so
$F(z)$ is a genuine polynomial with $\eta:=q-n-1\geq 0$,
\[
F(z)=F_l(z)=z^{n+\eta}B_1+~\ldots~+z^{1+\eta}B_n
\quad\quad\quad\quad \eta\geq 0.
\]
Then, the corresponding $p(n+\eta)\times m(n+\eta)~$ Hankel matrix
$\hat{\mathbf H}_{\mathbf\eta}$ (hat for left polynomial)
takes the form of
\begin{equation}\label{eq:ContrL}
\hat{\mathbf H}_{\mathbf\eta}=\left(\begin{smallmatrix}
~     &~     &B_n   &\cdots&B_1\\
~     &\Ddots&   ~  &\Ddots&   \\
B_n   &\cdots&B_1   & ~    & ~ \\
\vdots&\Ddots&~     & ~    & ~ \\
B_1&~ & ~    & ~    & ~    & ~~
\end{smallmatrix}\right).
\end{equation}
(iii) In the intermediate case, where $q\in[1,~n]$, $F(z)$ in
\eqref{eq:Poly} is a genuine Laurent polynomial and we shall
write it in the form of \eqref{eq:rational} (with $B_q=D$) as
\begin{equation}\label{eq:Laurent}
\begin{smallmatrix}
F(z)&=&F_l(z)+D+F_r(z)\\~\\
F_r(z)&=&z^{-1}B_{q+1}+z^{-2}B_{q+2}+~\ldots~+z^{q-n}B_n\\~\\
F_l(z)&=&z^{q-1}B_1+z^{q-2}B_2+~\ldots~+z^1B_{q-1}~.
\end{smallmatrix}
\quad\quad\quad\quad q\in[1,~n].
\end{equation}
The corresponding Hankel matrices for $F_r(z)$ and
$F_l(z)$, respectively are
\begin{equation}\label{eq:Hankels}
{\mathbf H}=\left(
\begin{smallmatrix}
B_{q+1}&B_{q+2}&\cdots&B_{n-1}&B_n\\
B_{q+2}&B_{q+3}&\cdots&B_n    & ~ \\
\vdots &   ~   &\Ddots&  ~    & ~ \\
B_{n-1}&B_n    &  ~   &~      & ~ \\
B_n    &~      &   ~  & ~     & ~
\end{smallmatrix}\right)
\quad\quad\quad\quad
\hat{\mathbf H}=\left(
\begin{smallmatrix}
B_{q-1}&B_{q-2}&\cdots&B_2&B_1\\
B_{q-2}&B_{q-3}&\cdots&B_1& ~ \\
\vdots &     ~ &\Ddots&~  & ~ \\
B_2    &B_1    &   ~  & ~ & ~ \\
B_1    &  ~    &   ~  & ~ & ~
\end{smallmatrix}\right).
\end{equation}

It is well known, see e.g. \cite[Theorem 4.5]{BGK}, that the McMillan
degree of $F(z)$ in \eqref{eq:rational} and in \eqref{eq:Laurent} is
equal to the sum of the ranks of $\hat{\mathbf H}$ and
${\mathbf H}$, the associated Hankel matrices.
Thus we can now state the main result of this subsection.
\vskip 0.2cm

\begin{Ob}\label{Ob:RankDeg}
Let us denote by $~d~$ the McMillan degree of the $p\times m$-valued
(possibly Laurent) polynomial $F(z)$ in \eqref{eq:Poly},
\[
F(z)=z^{q}\left(z^{-1}B_1+~\ldots~+z^{-n}B_n\right).
\]
Let ${\mathbf H}$, $\hat{\mathbf H}$ be the Hankel matrices
associated with $F_r(z)$, $F_l(z)$, respectively. Then,
\[
d=\left\{\begin{smallmatrix}
{\rm rank}({\mathbf H})&~&0\geq q&~&{\rm Eq.}\quad
\eqref{eq:ContrR}\\~\\
{\rm rank}(\hat{\mathbf H})&~&q\geq n+1&~&{\rm Eq.}\quad
\eqref{eq:ContrL}\\~\\
{\rm rank}({\mathbf H})+{\rm rank}(\hat{\mathbf H})&~&
q\in[1,~n]&~&{\rm Eq.}\quad\eqref{eq:Hankels}.
\end{smallmatrix}\right.
\]
\end{Ob}
\vskip 0.2cm

So far, we have used the Hankel matrix $\mathbf{H}$ to obtain the
McMillan degree of a given FIR $F(z)$. In Subsection
\ref{subsec:IsometricHankel} below we show that this $F(z)$
is in $\mathcal{U}$ if and only if $\mathbf{H^*}\mathbf{H}$ (or
$\mathbf{H}\mathbf{H^*}$) have a certain invariant subspace.

\subsection{Families of FIR systems in $\mathcal{U}$}
\label{subsec:families}

In this subsection we show how to produce, out of a given
$p\times m$-valued Laurent polynomial,
\begin{equation}\label{eq:OrigPoly}
F(z)=z^{q}\left(z^{-1}B_1+~\ldots~+z^{-n}B_n\right),
\end{equation}
a whole family of Laurent polynomials of ~{\em
various dimensions and powers}. Moreover, whenever the original
one is in ${\mathcal U}$, then so is all the resulting
family of polynomials.
\vskip 0.2cm

Clearly, when $F(z)$ in \eqref{eq:OrigPoly} is  in $\mathcal{U}$,
so are $F(z)U_m$ and $U_pF(z)$, where $U_m$ and $U_p$ are
arbitrary $m\times m$ and $p\times p$ constant unitary matrices.
\vskip 0.2cm

Note also that if for some value of the parameter $q$, $F(z)$ in
\eqref{eq:OrigPoly} is  in $\mathcal{U}$, then this is the case
for all $q$. Thus, without loss of generality, we can take all
polynomials to be causal, i.e. in \eqref{eq:OrigPoly} $0\geq q$.
\vskip 0.2cm

We here illustrate six versions of the newly generated polynomials.
Obviously, to further enrich the variety, they may be combined. 
\vskip 0.2cm

{\bf I.}\quad 
The reverse polynomial,
\[
\begin{matrix}
F_{\rm rev}(z)&:=&
z^{-1}B_n+z^{-2}B_{n-1}+~\ldots~+z^{1-n}B_2+z^{-n}B_1~.
\end{matrix}
\]
{\bf II.}\quad Preserving the McMillan degree.\\
For $-1\geq q$ and $j=1,~\ldots~,~1-q$, one can generate
$~jp\times jm$-valued polynomials. For example,
taking in \eqref{eq:OrigPoly} $n=4$ one may obtain:\\
{\bf a.}\quad For $q=-1$
\[
\begin{matrix}
F(z)&=&z^{-1}(z^{-1}B_1+z^{-2}B_2+z^{-2}B_3+z^{-4}B_4)&~&
p\times m{\rm-~valued}\\~\\
F(z)&=&z^{-1}\left(\begin{smallmatrix}0&B_1\\
B_1&B_2\end{smallmatrix}\right)+z^{-2}\left(\begin{smallmatrix}
B_2&B_3\\B_3&B_4\end{smallmatrix}\right)
+z^{-3}\left(\begin{smallmatrix}B_4&0\\ 0&0
\end{smallmatrix}\right)&~&2p\times 2m{\rm-~valued},
\end{matrix}
\]
{\bf b.}\quad For $q=-2$ 
\[
\begin{matrix}
F(z)&=&z^{-2}(z^{-1}B_1+z^{-2}B_2+z^{-2}B_3+z^{-4}B_4)&~&
p\times m{\rm-~valued}\\~\\
F(z)&=&z^{-1}\left(\begin{smallmatrix}0&0\\
0&B_1\end{smallmatrix}\right)+z^{-2}\left(\begin{smallmatrix}
B_1&B_2\\B_2&B_3\end{smallmatrix}\right)
+z^{-3}\left(\begin{smallmatrix}B_3&B_4\\ B_4&0
\end{smallmatrix}\right)&~&2p\times 2m{\rm-~valued}\\~\\
F(z)&=&z^{-1}\left(\begin{smallmatrix}0&0&B_1\\
0&B_1&B_2\\B_1&B_2&B_3\end{smallmatrix}\right)
+z^{-2}\left(\begin{smallmatrix}B_2&B_3&B_4\\
B_3&B_4&0\\B_4&0&0
\end{smallmatrix}\right)&~&3p\times 3m{\rm-~valued}.
\end{matrix}
\]
We shall see that both polynomials in ~~{\bf a}~~
share the same McMillan degree. A similar statement holds
for the three polynomials in ~{\bf b}.
\vskip 0.2cm

{\bf III.}\quad Doubling the powers.\\
For a pair of parameters $~a~$ integer
and $~\gamma~$ natural,
\begin{equation}\label{eq:z^a}
F(z)=z^a\left(z^{-1\cdot\gamma}B_1+z^{-2\cdot\gamma
}B_2+~\ldots~+z^{-n\cdot\gamma}B_n\right).
\end{equation}
In particular, for $a=0$, $F(z)$ in \eqref{eq:z^a} may be written
as $F(z^{\gamma})$. Rational functions within $\mathcal{U}$, of this
structure, suit what is known in signal processing as~
{\em filter banks}, see e.g. the books \cite{SN}, \cite{Vaid} and
the papers \cite{AJL3}, \cite{AJL5}, \cite{CZZM3} \cite{OTHN},
\cite{HTN}, \cite{RMW}, \cite{VHEK}, \cite{YaZh08}.
\vskip 0.2cm

In a similar way, one can obtain richer structures, e.g.
of the form
\begin{equation}\label{eq:jumps}
F(z)=z^{-1}B_1+z^{-2}B_2+z^{-3}B_3+z^{-9}B_4+z^{-10}B_5+z^{-11}B_6~.
\end{equation}
{\bf IV.}\quad Rectangular polynomials.
Let $\rho$ be a parameter so that $\frac{n}{\rho}$ is
natural\begin{footnote}{For given $n$ and $\rho$, one can
always find $\zeta\in [0, \rho-1]$ so that
$\frac{n+\zeta}{\rho}$ is natural. Then, the last part
of ${\mathbf B_{\frac{n}{\rho}}}$ in \eqref{eq:polyGamma}
is comprized of zeros.}\end{footnote},
\begin{equation}\label{eq:polyGamma}
F(z)=z^{a\rho}\left(z^{-1\cdot\gamma\rho}{\mathbf B_1}
+z^{-2\cdot\gamma\rho}{\mathbf B_2}
+~\ldots~+z^{-n\cdot\gamma\rho}{\mathbf B_{\frac{n}{\rho}}}
\right)
\end{equation}
{\bf a.}\quad $F(z)$ is a $~\rho p\times m$-valued
polynomial with coefficients,
\[
{\mathbf B_1}=\left(
\begin{smallmatrix}B_1\\ \vdots\\B_{\rho}\end{smallmatrix}\right)
\quad\quad
{\mathbf B_2}=
\left(\begin{smallmatrix}B_{\rho+1}\\ \vdots\\
B_{2\rho}\end{smallmatrix}\right)
\quad\quad
\cdots
\quad\quad
{\mathbf B_{\frac{n}{\rho}}}=
\left(\begin{smallmatrix}B_{n+1-\rho}\\
\vdots\\ B_{n}\end{smallmatrix}\right)
\]
{\bf b.}\quad $F(z)$ is a $~p\times\rho m$-valued
polynomial with coefficients,
\[
{\mathbf B_1}=\left(\begin{smallmatrix}B_1~\cdots~B_{\rho}\end{smallmatrix}\right)
\quad\quad
{\mathbf B_2}=
\left(\begin{smallmatrix}B_{\rho+1}~\cdots~B_{2\rho}\end{smallmatrix}\right)
\quad\quad
\cdots
\quad\quad
{\mathbf B_{\frac{n}{\rho}}}=
\left(\begin{smallmatrix}B_{n+1-\rho}~\cdots~B_{n}\end{smallmatrix}\right).
\]
{\bf V.}\quad Composition of polynomials.\\
Out of the pair of polynomials,
\begin{equation}\label{eq:PairPoly}
\begin{matrix}
F_b(z)&=&z^{-1}B_1+\ldots~z^{-n}B_n&&p_b\times m_b
\\~\\
F_c(z)&=&z^{-1}C_1+\ldots~z^{-l}C_l&&p_c
\times m_c
\end{matrix}\quad\quad\quad n\geq l,
\end{equation}
construct the following third polynomial
\[
F_d(z)=z^{-1}D_1+~\ldots~+z^{-n}D_n~,
\]
in three different forms.
\vskip 0.2cm

{\bf a.}\quad $F_d(z)$ is $(p_b+p_c)\times(m_b+m_c)$-valued
polynomial with coefficients,
\[
\begin{smallmatrix}
D_1=\left(\begin{smallmatrix}B_1&0\\~\\ 0& 
C_1\end{smallmatrix}\right)&&\cdots&&
D_l=\left(\begin{smallmatrix}B_l&0\\~\\ 
0&C_l\end{smallmatrix}\right)&&
D_{l+1}=\left(\begin{smallmatrix}B_{l+1}&0\\~\\ 
0&0_{p_c\times m_c}\end{smallmatrix}\right)&&
\cdots&&
D_n=\left(\begin{smallmatrix}B_n&0\\~\\
0&0_{p_c\times m_c}\end{smallmatrix}\right),
\end{smallmatrix}
\]
or
\[
\begin{smallmatrix}
D_1=\left(\begin{smallmatrix}0&B_1\\~\\ C_1&0
\end{smallmatrix}\right)
&&\cdots&&
D_l=\left(\begin{smallmatrix}0&B_l\\~\\
C_l&0\end{smallmatrix}\right)&&
D_{l+1}=\left(\begin{smallmatrix}0&B_{l+1}\\~\\
0_{p_c\times m_c}&0\end{smallmatrix}\right)
&&\cdots&&
D_n=\left(\begin{smallmatrix}0&B_n\\~\\
0_{p_c\times m_c}&0\end{smallmatrix}\right).
\end{smallmatrix}
\]
If both $F_b(z)$ and $F_c(z)$ are isometries on $\mathbb{T}$,
then so is $F_d(z)$.
\vskip 0.2cm

If both $F_b(z)$ and $F_c(z)$ are co-isometries on $\mathbb{T}$,
then so is $F_d(z)$.
\vskip 0.2cm

{\bf b.}\quad For $m_c\geq m_b$, $F_d(z)$ is
$~(p_b+p_c)\times m_c$-valued polynomial with coefficients,
\[
\begin{smallmatrix}
D_1
=\left(\begin{smallmatrix}\sqrt{\alpha}B_1~~0_{p_b\times(m_c-m_b)}
\\~\\ \sqrt{1-\alpha}C_1\end{smallmatrix}\right)&&
\cdots&&
D_l
=\left(\begin{smallmatrix}\sqrt{\alpha}B_l~~0{p_b\times(m_c-m_b)}
\\~\\ \sqrt{1-\alpha}C_L\end{smallmatrix}\right)\\~\\~\\
D_{l+1}
=\left(\begin{smallmatrix}\sqrt{\alpha}B_{l+1}~~0{p_b\times(m_c-m_b)}
\\~\\ 0_{p_c\times m_c}\end{smallmatrix}\right)&&
\cdots&&
D_n=\left(\begin{smallmatrix}\sqrt{\alpha}B_n~~0{p_b\times(m_c-m_b)}
\\~\\ 0_{p_c\times m_c}\end{smallmatrix}\right)
\end{smallmatrix}\quad\quad\alpha\in[0, 1]\quad{\rm parameter}.
\]
If both $F_b(z)$ and $F_c(z)$ are isometries on $\mathbb{T}$,
then so is $F_d(z)$
\vskip 0.2cm

{\bf c.}\quad For $p_b\geq p_c$, $F_d(z)$ is
$~p_b\times(m_b+m_c)$-valued polynomial with coefficients
\[
\begin{smallmatrix}
D_1=\left(\begin{smallmatrix}\begin{smallmatrix}
\sqrt{\alpha}B_1\\ 0_{(p_b-p_c)\times m_b}\end{smallmatrix}
\quad \sqrt{1-\alpha}C_1\end{smallmatrix}\right)&&
\cdots&&
D_l
=\left(\begin{smallmatrix}\begin{smallmatrix}
\sqrt{\alpha}B_l\\ 0_{(p_b-p_c)\times m_b}\end{smallmatrix}
\quad 
\sqrt{1-\alpha}C_l\end{smallmatrix}\right)\\~\\~\\
D_{l+1}
=\left(\begin{smallmatrix}\begin{smallmatrix}
\sqrt{\alpha}B_{l+1}\\ 0_{(p_b-p_c)\times m_b}\end{smallmatrix}
\quad 0_{p_b\times m_c}\end{smallmatrix}\right)&&
\cdots&&
D_n
=\left(\begin{smallmatrix}\begin{smallmatrix}
\sqrt{\alpha}B_n\\ 0_{(p_b-p_c)\times m_b}\end{smallmatrix}
\quad 0_{p_b\times m_c}\end{smallmatrix}\right)
\end{smallmatrix}
\quad\quad\alpha\in[0, 1]\quad{\rm parameter}.
\]
If both $F_b(z)$ and $F_c(z)$ are co-isometries on $\mathbb{T}$,
then so is $F_d(z)$
\vskip 0.2cm

{\bf VI.}\quad Product of the polynomials.
\[
\begin{matrix}
F_b(z)&=&z^{-1}B_1+\ldots~z^{-n}B_n&&p_b\times\rho
\\~\\
F_c(z)&=&z^{-1}C_1+\ldots~z^{-l}C_l&&\rho\times m_c
\end{matrix}
\]
i.e. $F_d(z)$ is the following $~p_b\times m_c$-valued polynomial
\[
F_b(z)F_c(z):=F_d(z)=z^{-1}\left(z^{-1}D_1+~\ldots~+z^{-(n+l-1)}D_{n+l-1}\right),
\]
where a straightforward computation yields that the coefficients
$D_1~,~\ldots~,~D_{n+l-1}$ are given by,
\begin{equation}\label{eq:CoeffD}
{\mathbf D}_{1, n+l-1}=\mathbf{H}_l{\mathbf T}_{n+l,~\rho}
{\mathbf C}_{1, l}
\end{equation}
with ${\mathbf C}_{1, l}$ and ${\mathbf D}_{1, n+l-1}$ in the spirit
of \eqref{eq:B_eta}, $\mathbf{H}_l$ is the $(n+l)p\times(n+l)m$ Hankel
matrix as in \eqref{eq:ContrR} (with $\eta=l$) and 
\begin{equation}\label{eq:T}
{\mathbf T}_{n+l,~\rho}:=\left(\begin{smallmatrix}
~&~&~&~&~&I_{\rho}\\
~&~&~&~&I_{\rho}&~\\
~&~&\Ddots&~&~\\
~&I_{\rho}&~&~&~&~\\
I_{\rho}&~&~&~&~&~\end{smallmatrix}\right).
\end{equation}
In the Appendix below we show how each of the new polynomials, in the
above items {\bf I} through {\bf VI}, is constructed.
A key tool there will be the associated Hankel matrix.
\vskip 0.2cm

Furthermore, as already pointed out, we shall see that if the original
polynomial $F(z)$ in \eqref{eq:OrigPoly} is in $\mathcal{U}$, so are
all the resulting polynomials.
\vskip 0.2cm

To this end, in Subsections \ref{sec:RealizeMatrix},
\ref{sec:Blaschke} and \ref{subsec:IsometricHankel} respectively,
we present three characterizations of Laurent polynomials in $\mathcal{U}$.

\section{minimal realization of (co)-isometric FIR}
\setcounter{equation}{0}
\label{sec:IsoRealization}

\subsection{Characterization through realization matrices of 
Schur stable systems}
\label{sec:RealizeMatrix}

In this subsection we characterize, through a corresponding realization
matrix, rational functions, analytic outside the open unit disk (Schur
stable), which are in $\mathcal{U}$.
\vskip 0.2cm

There are several variants of characterizations of $F\in\mathcal{U}$
through its minimal realization matrix $R$. The square case ($m=p$)
was addressed in \cite[Lemma 2 \& Theorem 3]{GVKDM}. Another version
of it appeared in \cite[Theorem 3.1]{AlGo1} and
\cite[Theorem 2.1]{AlGo2}.
Below we cite and adapted form of \cite[Theorem 14.5.1]{Vaid}. A more
general case was studied in \cite[Theorem 4.5]{AlRak2}. In fact, in
\cite{AlGo1}, \cite{AlGo2}, \cite{AlRak2} and \cite{GVKDM}
they address the ~{\em indefinite}~ inner product case where
$\left(F(z)\right)^*J_pF(z)=J_m$ with $J_p, J_m$ signature matrices, i.e. diagonal
matrices satisfying $J_p^2=I_p$ and $J_m^2=I_m$.

\begin{Tm}\label{Tm:RealizationU}
Let $F(z)$ be a $p\times m$-valued rational function whose poles
are within the open unit disk and let $R$ be a corresponding
$(\nu+p)\times(\nu+m)$ minimal realization matrix \eqref{eq:R}
\[
R:={\footnotesize
\left(\begin{array}{c|c}A&B \\ \hline C&D
\end{array}\right)}.
\]
I. Assume that $p\geq m$.
$F(z)$ is 
in $\mathcal{U}$ if and only if,
\begin{equation}\label{eq:SteinRiso}
R^*{\rm diag}\{I_{\nu}\quad I_m\}R={\rm diag}\{I_{\nu}\quad I_m\}.
\end{equation}
II. Assume that $m\geq p$.
$F(z)$ is 
in $\mathcal{U}$ if and only if,
\begin{equation}\label{eq:SteinRcoIso}
R{\rm diag}\{I_{\nu}\quad I_p\}R^*={\rm diag}\{I_{\nu}\quad I_p\}.
\end{equation}
\end{Tm}
\vskip 0.2cm

We now recall the notion of Controllability and Observability
Gramians, see e.g. \cite[Eqns (12.8.17), (12.8.37)]{HSK}.
\vskip 0.2cm

Assuming that the spectrum of $A$ (the upper left block in $R$ in
\eqref{eq:R}) is within the open unit disk (Schur stable), 
$W_{\rm cont}$, $W_{\rm obs}$, the associated $\nu\times\nu$
Controllability and Observability Gramians, respectively, are
given by the solution to the respective Stein equations
\begin{equation}\label{eq:gramian}
W_{\rm cont}-AW_{\rm cont}A^*=BB^*
\quad\quad\quad\quad
W_{\rm obs}-A^*W_{\rm obs}A=C^*C.
\end{equation}
We can now state the following whose proof is given in
\cite{AJL4}.

\begin{La}\label{CorollaryGramians}
Let $F(z)$ be a $p\times m$-valued rational function whose poles
are within the open unit disk and denote by
$W_{\rm cont}$, $W_{\rm obs}$ the associated controllability and
observability gramians, respectively.
\vskip 0.2cm

Assume that $F(z)$ is in $\mathcal{U}$.
\vskip 0.2cm

I. If $p\geq m$, $F(z)$ admits a state space realization $R$ in
\eqref{eq:SteinRiso} so that
\[
(I_{\nu}-W_{\rm cont})\quad{\rm positive~~semidefinite}
\quad\quad\quad\quad\quad W_{\rm obs}=I_{\nu}~.
\]
II. If $m\geq p$, $F(z)$ admits a state space realization $R$ in
\eqref{eq:SteinRcoIso} so that
\[
W_{\rm cont}=I_{\nu}\quad\quad\quad\quad\quad
(I_{\nu}-W_{\rm obs})\quad{\rm positive~~semidefinite}.
\]
III. If $p=m$, $F(z)$ admits a state space realization $R$
in \eqref{eq:SteinRiso}, \eqref{eq:SteinRcoIso}
so that
\[
W_{\rm cont}=I_{\nu}\quad\quad\quad\quad\quad W_{\rm obs}=I_{\nu}~.
\]
\end{La}
\vskip 0.2cm

This result will be used in the proof of Proposition \ref{Pn:UnitaryF}
below. 
\vskip 0.2cm

We conclude this subsection by illustrating the results of Observation
\ref{Ob:RankDeg}, Theorem \ref{Tm:RealizationU} and of Lemma
\ref{CorollaryGramians}
\vskip 0.2cm

\begin{Ex}\label{Ex:Polynomials}
{\rm 
{\bf I.}\quad
Consider the $2\times 2$-valued polynomial $F(x)$ in
\eqref{eq:poly} with $n=3$ i.e.
\[
F(z)=z^q(z^{-1}B_1+z^{-2}B_2+z^{-3}B_3),
\]
where
\begin{equation}\label{ExampA1}
B_1=
\scriptstyle{\frac{1}{5}}
\left(\begin{smallmatrix}2&2\\ 2&2\end{smallmatrix}\right)
\quad\quad\quad\quad
B_2=\scriptstyle{\frac{1}{5}}
\left(\begin{smallmatrix}~~~0&3\\-3&0\end{smallmatrix}\right)
\quad\quad\quad\quad
B_3=\scriptstyle{\frac{1}{5}}
\left(\begin{smallmatrix}~~~2&-2\\-2&~~~2\end{smallmatrix}\right),
\end{equation}
and $q$ is a parameter assuming the values 2~ and ~1.
\vskip 0.2cm

{\bf (i)}\quad For $q=2$, $F(z)$ takes the form,
\[
F(z)=zB_1+B_2+z^{-1}B_3.
\]
Namely, as in \eqref{eq:Laurent}, and here \eqref{eq:Hankels}
takes the form of ${\mathbf H}=B_3$ and $\hat{\mathbf H}=B_1$
so the McMillan degree is
$d={\rm rank}(\hat{\mathbf H})+{\rm rank}({\mathbf H})=2$.
\vskip 0.2cm

A minimal realization (i.e. of dimension 2) is
\[
F(z)=zv_1v_1^*+B_2+\frac{1}{z}v_3v_3^*
\quad\quad{\rm with}\quad\quad
v_1=\scriptstyle{\frac{\sqrt 2}{\sqrt 5}}
\left(\begin{smallmatrix}1\\~\\ 1\end{smallmatrix}\right)
\quad
v_3=\scriptstyle{\frac{\sqrt 2}{\sqrt 5}}
\left(\begin{smallmatrix}~~~1\\~\\-1\end{smallmatrix}\right).
\]

{\bf (ii)}\quad For $q=1$, $F(z)$ is causal and takes the form,
\[
F(z)=B_1+z^{-1}B_2+z^{-2}B_3.
\]
A straightforward dimension 4 realization of the form
\eqref{eq:StraightforwardRealiz} is
\[
\tilde{R}_r={\footnotesize\left(\begin{array}{cc|c}
0&I_2&B_2\\ 0&0&B_3\\ \hline
I_2&0&B_1\end{array}\right)}.
\]
A closer scrutiny reveals that 
here \eqref{eq:Hankels}
takes the form of $\hat{\mathbf H}$ empty, and
\[
{\mathbf H}=
\left(\begin{smallmatrix}B_2&B_3\\B_3&0\end{smallmatrix}\right)=
\scriptstyle{\frac{1}{5}}\left(\begin{smallmatrix}
~~~0&~~~3&~~~2&-2\\
-3&~~0&-2&~~2\\
~~~2&-2&~~~0&~~~0\\
-2&~~~2&~~~0&~~~0\end{smallmatrix}\right).
\]
Thus the McMillan degree is
$d={\rm rank}(\hat{\mathbf H})+{\rm rank}({\mathbf H})=0+2=2$.
\vskip 0.2cm

Indeed, a minimal realization is,
\[
R_r=\scriptstyle{\frac{1}{5}}
{\footnotesize\left(\begin{array}{rr|rr}
-2&-2&1&-4\\
~~~2&~~~2&4&-1\\ \hline
-4&~~~1&2&~~~2\\~~~1&-4&2&~~~2\end{array}\right)}.
\]
Recall that from both parts of Theorem \ref{Tm:RealizationU}
it follows that 
\mbox{$R_rR_r^*={\rm diag}\{I_2~,~I_2\}=R_r^*R_r$}.
\vskip 0.2cm

{\bf II.}\quad 
Consider the $1\times 2$-valued polynomial
\begin{equation}\label{eq:CoIsoEx}
F(z)=z^q\left(z^{-1}
\left(0\quad{\scriptstyle-\frac{3}{5}}\right)+
z^{-2}\left({\scriptstyle\frac{4}{5}}\quad 0\right)\right),
\end{equation}
where the parameter $q$ assume the values 0 and 1.
\vskip 0.2cm

For $q=0$, $F(z)$ is strictly causal. The associated Hankel matrix
${\mathbf H}$ in \eqref{eq:ContrR} is given by,
\begin{equation}\label{eq:ExampHr}
{\mathbf H}=\scriptstyle{\frac{1}{5}}\left(
\begin{smallmatrix}
0&-3&~4&0\\
4&~~~0&~0&0\end{smallmatrix}\right)
\end{equation}
and the Hankel singular values are 1 and $0.8$.
A minimal realization is
\[
R={\footnotesize\left(\begin{array}{cc|cr}
0&\frac{4}{5}&0&-\frac{3}{5}\\
0&0&1&~~~0\\ \hline
1&0&0&0
\end{array}\right)}.
\]
Indeed, from part II of Theorem \ref{Tm:RealizationU} it follows
that $R~$ is co-isometry, i.e.  \mbox{$RR^*={\rm diag}\{I_2~,~1\}$}
and from Lemma \ref{CorollaryGramians}, we have that the
Observability Gramian is
\mbox{$W_{\rm obs}={\rm diag}\{1\quad\frac{16}{25}\}$}.
\vskip 0.2cm

Substituting in \eqref{eq:CoIsoEx} $q=1$ yields the causal
polynomial,
\[
F(z)=\left(0\quad{\scriptstyle-\frac{3}{5}}\right)+
z^{-1}\left({\scriptstyle\frac{4}{5}}\quad 0\right).
\]
A corresponding minimal realization is
\[
R={\footnotesize\left(\begin{array}{c|rr}
0&1&0\\ \hline\frac{4}{5}&0&-\frac{3}{5}
\end{array}\right)}.
\]
Indeed, $RR^*={\rm diag}\{1\quad 1\}$ and
$~W_{\rm obs}=\frac{16}{25}~$.
}
\qed
\end{Ex}
\vskip 0.2cm

In this subsection we have assumed that the system is causal
and in particular Schur stable. In the sequel, we remove this
restriction.

\subsection{A characterization through the Blaschke-Potapov product}
\label{sec:Blaschke}

We first need some preliminaries. We shall use the fact that
a $k\times k$ rank one orthogonal projection,
\[
P^*=P=P^2
\quad\quad\quad\quad{\rm rank}(P)=1,
\]
can always be written as
\begin{equation}\label{eq:DefProjVec}
P=vv^*\quad\quad v^*v=1\quad\quad v\in\C^k.
\end{equation}
Using \eqref{eq:DefProjVec}, a rank $k-1$ orthogonal projection $~Q~$ i.e.
\[
Q^*=Q^2=Q\quad\quad\quad\quad\quad\quad{\rm rank}(Q)=k-1,
\]
can always be written as
\begin{equation}\label{eq:Qprojection}
Q:=I_k-vv^*.
\end{equation}
We can now cite the classical Blaschke-Potapov product
result (as appeared in \cite{AJL4}). 

\begin{La}\label{Tm:U(z)}
A $~p\times m$-valued rational function $F(z)$, of McMillan
degree $~d$, is in $~{\mathcal U}$, \eqref{eq:DefU}, if and
only if
\begin{equation}\label{eq:UnitaryRational1}
\begin{matrix}
p\geq m&F(z)=\left(\prod\limits_{j=1}^d\left(I_p+
\left(\frac{1-{\alpha}^*_jz}{z-\alpha_j}-1\right)v_jv_j^*\right)
\right)U_{\rm iso}\\~\\
m\geq p&F(z)=U_{\rm coiso}\left(\prod\limits_{j=1}^d
\left(I_m+\left(\frac{1-{\alpha}^*_jz}{z-\alpha_j}-1\right)v_jv_j^*
\right)\right)
\end{matrix}\quad
\begin{smallmatrix}
v_j\in\C^p&
v_j^*v_j=1
\\~\\
U_{\rm iso}\in\C^{p\times m}&
U_{\rm iso}^*U_{\rm iso}=I_m
\\~\\
\quad\alpha_j\in\{\infty\cup\C\}\smallsetminus\mathbb{T}&~\\~\\
v_j\in\C^m&
v_j^*v_j=1
\\~\\
U_{\rm coiso}\in\C^{p\times m}&U_{\rm coiso}U_{\rm coiso}^*=I_p~.
\end{smallmatrix}
\end{equation}
\end{La}
\vskip 0.2cm

Recall $\prod\limits_{j=1}^0:=I$
\vskip 0.2cm

In Theorem \ref{Tm:UnitaryFIR} below,
we specialize this result to the FIR framework.
To this end, we recall that
\begin{equation}\label{eq:InvBla}
\left(I+\left(\frac{1-{\alpha}^*z}{z-\alpha}-1\right)vv^*\right)^{-1}
=\left(I+\left(\frac{z-\alpha}{1-{\alpha}^*z}-1\right)vv^*\right)
\quad
\begin{smallmatrix}
v^*v=1\\~\\
\alpha\in\{\infty\cup\C\}\smallsetminus\mathbb{T}.
\end{smallmatrix}
\end{equation}
We can now formulate the main result of this subsection.

\begin{Tm}\label{Tm:UnitaryFIR}
Let $F(z)$ be a $p\times m$-valued Finite Impulse Response function
of McMillan degree $~d$.
\vskip 0.2cm

I. Assuming $p\geq m$, $F(z)$ is isometrically in ${\mathcal U}$
if and only if 
it can be written in one, and hence in each, of these three
equivalent forms,
\begin{equation}\label{eq:FirIsoRational}
F(z)=\left\{\begin{smallmatrix}
\prod\limits_{j=1}^{\gamma}\left(I_p+(z-1)v_jv_j^*\right)
\prod\limits_{j=\gamma+1}^d\left(I_p+(\frac{1}{z}-1)v_jv_j^*
\right)U_{\rm iso}\\
\prod\limits_{j=1}^{\gamma}\left(I_p+(z-1)v_jv_j^*\right)
\left(\prod\limits_{j=d}^{\gamma+1}\left(I_p+(z-1)v_jv_j^*\right)
\right)^{-1}U_{\rm iso}
\\
\left(\prod\limits_{j=\gamma}^1\left(I_p+(\frac{1}{z}-1)v_jv_j^*
\right)\right)^{-1}
\prod\limits_{j=\gamma+1}^d\left(I_p+(\frac{1}{z}-1)v_jv_j^*
\right)U_{\rm iso}\end{smallmatrix}\right.
\quad\quad\begin{smallmatrix}
U_{\rm iso}^*U_{\rm iso}=I_m\\~\\
v_j\in\C^p\\~\\
v_j^*v_j=1\\~\\
\gamma\in[0,~d].
\end{smallmatrix}
\end{equation}
II. For $m\geq p$, $F(z)$ is co-isometrically in ${\mathcal U}$
if and only if it can be written in one, and hence in each, of
these three equivalent forms,
\begin{equation}\label{eq:FirCoIsoRational}
F(z)=\left\{\begin{smallmatrix}
U_{\rm coiso}\prod\limits_{j=\gamma+1}^d
\left(I_m+(\frac{1}{z}-1)v_jv_j^*\right)
\prod\limits_{j=1}^{\gamma}\left(I_m+(z-1)v_jv_j^*\right)
\\
U_{\rm coiso}\left(\prod\limits_{j=d}^{\gamma+1}\left(
I_m+(z-1)v_jv_j^*\right)\right)^{-1}
\prod\limits_{j=1}^{\gamma}\left(I_m+(z-1)v_jv_j^*\right)
\\
U_{\rm coiso}\prod\limits_{j=\gamma+1}^d\left(
I_m+(\frac{1}{z}-1)v_jv_j^*\right)
\left(\prod\limits_{j=\gamma}^1\left(
I_m+(\frac{1}{z}-1)v_jv_j^*\right)\right)^{-1}
\end{smallmatrix}\right.
\quad\begin{smallmatrix}
U_{\rm coiso}U_{\rm coiso}^*=I_p\\~\\
v_j\in\C^m\\~\\
v_j^*v_j=1\\~\\
\gamma\in[0,~d].
\end{smallmatrix}
\end{equation}
In particular, the function $F(z)$ in \eqref{eq:FirIsoRational},
\eqref{eq:FirCoIsoRational} is causal for $\gamma=0$ and
anti-causal for $\gamma=d$.
\end{Tm}
\vskip 0.2cm

The first version in both \eqref{eq:FirIsoRational} and
\eqref{eq:FirCoIsoRational} is a re-writing of
\eqref{eq:UnitaryRational1} with
\mbox{$\alpha_1=~\ldots~=\alpha_{\gamma}=\infty$} and
\mbox{$\alpha_{\gamma+1}=~\ldots~=\alpha_d=0$}.
Using \eqref{eq:InvBla}, the other two versions follow
from the first.
\vskip 0.2cm

Note that it is only for the causal case, $\gamma=0$, or the anti-causal
case, $\gamma=d$, that in \eqref{eq:FirIsoRational} a product of the form,
\begin{equation}\label{eq:LaurentExample}
\prod\limits_{j=1}^{\gamma}
\left(I_p+(z-1)v_jv_j^*\right)
\prod\limits_{j=\gamma+1}^d
\left(I_p+(\frac{1}{z}-1)v_jv_j^*\right)U_{\rm iso}
\quad\quad\begin{smallmatrix}
U_{\rm iso}^*U_{\rm iso}=I_m\\~\\
v_j\in\C^p\\~\\
v_j^*v_j=1\\~\\
\gamma\in[0,~d],
\end{smallmatrix}
\end{equation}
is of McMillan degree $~d~$ for ~{\em any choice}~ of the 
projections $v_1~,~\ldots~,~v_d$. For \mbox{$\gamma\in[1,~d-1]$}, 
the McMillan degree of the expression in \eqref{eq:LaurentExample}
is ~{\em at most}~ $d$. For example, repetitive use of
\eqref{eq:InvBla} reveals that
substituting in \eqref{eq:LaurentExample}
\[
v_{j}v_j^*=v_{d+1-j}v_{d+1-j}^*\quad\quad j=1,~\ldots~,\gamma
\quad\quad\quad d=2\gamma\quad\quad\gamma\quad{\rm natural},
\]
yields the {\em zero}~ degree function $~I_p$.
\vskip 0.2cm

Theorem \ref{Tm:UnitaryFIR} asserts that whenever $F\in\mathcal{U}$
is of McMillan degree $d$, {\em there exist}~ rank one orthogonal
projections $~v_1v_1^*~,~\ldots~,~v_dv_d^*~$ satisfying
\eqref{eq:FirIsoRational} and \eqref{eq:FirCoIsoRational}.
\vskip 0.2cm

Recall that employing optimization for design of FIRs is fairly
common, see e.g. \cite{CZZM3}, \cite{Dav}, \cite{HTN}, \cite{TV}
and \cite{VHEK} and for FIRs within $\mathcal{U}$ see e.g.
\cite{RMW}. This motivates a ~{\em convex}~ description
of large families of FIRs. To this end, we next specialize
\cite[Observation 4.3]{AJL4}.
\vskip 0.2cm

\begin{Cy}\label{FIRangles}
All $p\times m$-valued FIR rational functions of
McMillan degree $~d~$ in $\mathcal{U}$ may be parametrized by,
\begin{equation}\label{eq:FIRangle}
\begin{matrix}
(d+1)\cdot [0,~2\pi)^{(2p-m-1)(m+d)+d(m-1)+m
}
&~&p\geq m\\~\\
(d+1)\cdot [0,~2\pi)^{
(2m-p-1)(p+d)+d(p-1)+p}
&~&m\geq p.
\end{matrix}
\end{equation}
The causal subset (here, =lossless subset), may be parametrized by,
\begin{equation}\label{eq:angleCausal}
\begin{matrix}
[0,~2\pi)^{(2p-m-1)(m+d)+d(m-1)+m}
&~&p\geq m\\~\\
[0,~2\pi)^{(2m-p-1)(p+d)+d(p-1)+p}
&~&m\geq p.
\end{matrix}
\end{equation}
\end{Cy}

There are a few parametrizations of all FIRs in $\mathcal{U}$,
see e.g. \cite[Propri\'{e}t\'{e} 41]{Icart}, \cite[Section
14.4]{Vaid} and the more detailed study in \cite{GNS}.
The above choice is advantageous as the set of parameters in
\eqref{eq:angleCausal}
is not only ~{\em convex}, but in fact a ~{\em polytope}. This
corresponds to having in \eqref{eq:FirIsoRational} or in
\eqref{eq:FirCoIsoRational}, $\gamma=0$. As the integral parameter 
$\gamma$ attains values in $[0,~d]$, the parameter set
in \eqref{eq:FIRangle} are $d+1$ copies of this ~{\em polytope}.
This is in particular convenient if one
wishes to:\\
(i) Design FIRs within $\mathcal{U}$ through optimization,
see e.g. \cite{RMW}.\\
(ii) Iteratively apply para-unitary similarity, see e.g.
\cite[Section 3.3]{Icart}, \cite{MWBCR}, \cite{SDeLID}.
In signal processing literature, this is associated with
with {\em channel equalization}~ and in communications
literature with ~{\em decorrelation of signals}~ or\\
(iii) Iteratively apply Q-R factorization in the
framework of communications, see e.g. \cite{CB1}, \cite{CB2}.
\vskip 0.2cm

We now establish a connection between 
Blaschke-Potapov description of $F(z)$ in \eqref{eq:FirIsoRational},
\eqref{eq:FirCoIsoRational} (with $\gamma=0$) and $B_1,~\ldots~,~B_n$
the coefficients of the polynomial
$F(z)$ in \eqref{eq:Poly} with $q=1$, i.e.
\begin{equation}\label{eq:q=1}
F(z)=B_1+z^{-1}B_2+~\ldots~z^{1-n}Z_n~.
\end{equation}

\begin{Ob}\label{Ob:Coefficients}
Assume that the (causal) $p\times m$-valued polynomial $F(z)$ in
\eqref{eq:q=1} is in $\mathcal{U}$ and of McMillan degree $~d$.
\vskip 0.2cm

Then if $m=p$ up to multiplication by a unitary matrix from the
left or from the right, the square coefficient matrices
$B_1~,~\ldots~,~B_n$, are given by,
\[
\begin{smallmatrix}
B_1&=&Q_1\cdots Q_d\\
B_2&=&\sum\limits_{j=1}^dQ_1\cdots Q_{j-1}v_jv_j^*Q_{j+1}\cdots Q_d\\
B_3&=&\sum\limits_{k=j+1}^d\sum\limits_{j=1}^{d-1}
Q_1\cdots Q_{j-1}v_jv_j^*Q_{j+1}\cdots Q_{k-1}v_kv_k^*Q_{k+1}\cdots Q_d\\
B_4&=&\sum\limits_{q=k+1}^d\sum\limits_{k=j+1}^{d-1}
\sum\limits_{j=1}^{d-2}Q_1\cdots Q_{j-1}v_jv_j^*Q_{j+1}\cdots
Q_{k-1}v_kv_k^*Q_{k+1}\cdots Q_{q-1}v_qv_q^*Q_{q+1}\cdots Q_d\\
\vdots\\~\\
B_n&=&v_1v_1^*\cdots~v_dv_d^*
\end{smallmatrix}
\]
where $~v_1v_1^*,~\ldots~,~v_dv_d^*$, are rank one orthogonal projections
\eqref{eq:DefProjVec} and \mbox{$Q_j:=I-v_jv_j^*$},
\mbox{$j=1,~\ldots~,~d$}, see \eqref{eq:Qprojection}.
\vskip 0.2cm

If $p\geq m$, the above $p\times p$ coefficient matrices
$B_1~,~\ldots~,~B_n$ are multiplied from the right by a
$~p\times m$ isometry $U_{\rm iso}$,
$(U_{\rm iso}^*U_{\rm iso}=I_m)$.
\vskip 0.2cm

If $m\geq p$, the above $m\times m$ coefficient matrices
$B_1~,~\ldots~,~B_n$ are multiplied from the left by a
$~p\times m$ co-isometry $U_{\rm coiso}$,
$(U_{\rm coiso}U_{\rm coiso}^*=I_p)$.
\end{Ob}
\vskip 0.2cm

In particular this implies that
\[
F(1)=B_1~+\ldots~+B_n=\left\{\begin{smallmatrix}
U_{\rm iso}&~&F~{\rm isometry}\\~\\
U_{\rm coiso}&~&F~{\rm co-isometry}
\end{smallmatrix}\right.
\]
\vskip 0.2cm

We conclude this section by pointing out that
so far we have focused on characterizations of Finite Impulse
Response functions in $\mathcal{U}$ through their ~{\em
minimal realization}.~ In the sequel
this restriction is removed.

\section{Hankel matrices - revisited}
\label{sec:Hankelrevisited}
\setcounter{equation}{0}

\subsection{Isometric FIR - a Hankel matrix characterization}
\label{subsec:IsometricHankel}

\vskip 0.2cm

Let $F(z)$, \eqref{eq:poly}, be a Finite Impulse Response (possibly
rectangular) rational function,
\[
F(z)=z^q\left(z^{-1}B_1+~\ldots~,~z^{-n}B_n\right)
\quad\quad\quad\quad q\quad{\rm parameter}.
\]
Clearly, the above $F(z)$ is in $\mathcal{U}$, if and only if
$F_o(z)$ in \eqref{eq:HatFr},
\[
F(z)_{|_{q=0}}=F_o(z):=z^{-1}B_1+~\ldots~,~z^{-n}B_n~,
\]
is in $\mathcal{U}$.
\vskip 0.2cm

Thus, in the sequel we find it convenient to focus on $F_o(z)$.
Subsequently, substituting in \eqref{eq:B_eta} and in \eqref{eq:ContrR} 
$\eta=0$ one obtains the matrices
\[
{\mathbf B}_0=\left(\begin{smallmatrix}
B_1\\
\vdots\\~\\
B_n
\end{smallmatrix}\right)
\quad\quad\quad\quad
{\mathbf H}_0=\left(\begin{smallmatrix}
B_1    &\cdots &B_n\\
\vdots &\Ddots &  ~\\
B_n    &~      & ~
\end{smallmatrix}\right).
\]
We shall also find it convenient to write $F_o(z)$ in
\eqref{eq:HatFr} in each of the two following forms
\begin{equation}\label{eq:F=ZB}
F_o(z)={\mathbf Z}{\mathbf B}_0
=\hat{\mathbf B}\hat{\mathbf Z},
\end{equation}
where ${\mathbf B}_0$ is as before and
\begin{equation}\label{eq:HatHatB}
\hat{\mathbf B}:=(B_1~,~\cdots~,~B_n),
\end{equation}
and 
\begin{equation}\label{eq:Z}
\begin{smallmatrix}
{\mathbf Z}&:=&(z^{-1}~\cdots~z^{-n})\otimes{I}_p&=&
(z^{-1}I_p~\cdots~z^{-n}I_p)\\~\\
\hat{\mathbf Z}&:=&\left(\begin{smallmatrix}z^{-1}\\\vdots\\~
\\z^{-n}\end{smallmatrix}\right)\otimes{I}_m&=&
\left(\begin{smallmatrix}z^{-1}I_m\\ 
\vdots \\~\\z^{-n}I_m\end{smallmatrix}\right).
\end{smallmatrix}
\end{equation}
where $\otimes$ denotes the usual Kronecker (=tensor) product,
see e.g. \cite[Section 4.2]{HJ2}.
\vskip 0.2cm

We now state the main result of this subsection.
\vskip 0.2cm

\begin{Tm}\label{Ob:UnitaryF}
Let $F_0(z)$ be a $p\times m$ polynomial in \eqref{eq:HatFr}
\[
F_0(z)=z^{-1}B_1+~\ldots~+z^{-n}B_n,
\]
and let ${\mathbf H}_0$ be the corresponding Hankel matrix.
\vskip 0.2cm

The polynomial $F_o(z)$ is in $\mathcal{U}$ if and only if 
\begin{equation}\label{eq:H^*H}
\begin{matrix}
\left(I_{nm}-{\mathbf H}_0^*{\mathbf H}_0\right)
\left(\begin{smallmatrix}
I_m\\0_{m(n-1)\times m}\end{smallmatrix}\right)&=&0_{pn\times m}&&&p\geq m\\~\\
\left(I_p\quad 0_{p\times(n-1)p}\right)
\left(I_{np}-{\mathbf H}_0{\mathbf H}_0^*\right)
&=&0_{p\times np}&&&m\geq p.
\end{matrix}
\end{equation}
\end{Tm}

{\bf Proof:}\quad A straightforward substitution of
$F_o(z)$ in the definition of $\mathcal{U}$ \eqref{eq:DefU}
yields the following characterization,
\begin{equation}\label{eq:ClassicalCharacter}
\begin{matrix}
\left(\begin{smallmatrix}
\sum\limits_{j=1}^nB_j^*B_j\\
\sum\limits_{j=1}^{n-1}B_{1+j}^*B_j\\
\vdots\\
\sum\limits_{j=1}^2B_{n-2+j}^*B_j\\
B_n^*B_1
\end{smallmatrix}\right)
&=&\left(\begin{smallmatrix}
I_m\\ 0_{m(n-1)\times m}
\end{smallmatrix}\right)&~&~&
p\geq m\\~\\~\\
\left(
\begin{smallmatrix}
\sum\limits_{j=1}^nB_jB_j^*~,
&
\sum\limits_{j=1}^{n-1}B_{1+j}B_j^*~,&
\cdots~,&\sum\limits_{j=1}^2B_{n-2+j}B_j^*~,&
B_nB_1^*\end{smallmatrix}\right)&=&
\left(\begin{smallmatrix}I_p~,&
0_{p\times(n-1)p}
\end{smallmatrix}\right)&~&~&
m\geq p.
\end{matrix}
\end{equation}
For the square case, see e.g. \cite[Propri\'{e}t\'{e} 37]{Icart}.
\vskip 0.2cm

Using the matrices $~{\mathbf H}_0~$ and ${\mathbf B}_0$ 
this may be equivalently written as,
\begin{equation}\label{eq:H^*B}
\begin{matrix}
{\mathbf H}_0^*{\mathbf B}_0=\left(\begin{smallmatrix}
I_m\\ 0_{m(n-1)\times m}
\end{smallmatrix}\right)&~&p\geq m,\\~\\
\hat{\mathbf B}
{\mathbf H}_{\mathbf 0}^*=(I_p\quad
0_{p\times(n-1)p}
)&~&m\geq p.
\end{matrix}
\end{equation}
where $\hat{\mathbf B}$ is as in \eqref{eq:HatHatB}.
\vskip 0.2cm

Now, since
\[
\begin{matrix}
{\mathbf B}_0={\mathbf H}_0
\left(\begin{smallmatrix}I_m\\ 0_{m(n-1)\times m}
\end{smallmatrix}\right)&&&{\rm and}&&&
\hat{\mathbf B}=(I_p\quad 0_{p\times(n-1)p}){\mathbf H}_0
\end{matrix}
\]
the relation in \eqref{eq:H^*B} is equivalent to
\[
\begin{matrix}
{\mathbf H}_0^*{\mathbf H}_0\left(\begin{smallmatrix}
I_m\\0_{m(n-1)\times m}\end{smallmatrix}\right)&=&
\left(\begin{smallmatrix}
I_m\\0_{m(n-1)\times m}\end{smallmatrix}\right)\\~\\
\left(I_p\quad 0_{p\times(n-1)p}\right){\mathbf H}_0{\mathbf H}_0^*
&=&\left(I_p\quad 0_{p\times(n-1)p}\right),
\end{matrix}
\]
which in turn can be written as \eqref{eq:H^*H}.
\qed
\vskip 0.2cm

An alternative proof of the same result is given in
\cite[Theorem 5.2]{AJL4}.
\vskip 0.2cm

Theorem \ref{Ob:UnitaryF} offers a characterization of a Laurent
polynomial $~F(z)~$ in $~\mathcal{U}$ through the existence of
a certain invariant subspace of ${\mathbf H}_0^*{\mathbf H}_0$
(or ${\mathbf H}_0{\mathbf H}_0^*$), see \eqref{eq:H^*H}.
A closer scrutiny reveals that whenever $~F(z)~$
is in $~\mathcal{U}$, more can be said about these matrices.
\vskip 0.2cm

\begin{Pn}\label{Pn:UnitaryF}
Let $F_o(z)$ be a ~$p\times m$-valued polynomial in \eqref{eq:HatFr}
\[
F_o(z)=z^{-1}B_1+~\ldots~+z^{-n}B_n,
\]
and let ${\mathbf H}_0$ be the corresponding Hankel matrix.
\vskip 0.2cm

Assume that $F(z)$ is in $\mathcal{U}$.
\vskip 0.2cm

I. For $p\geq m$, 
\begin{equation}\label{eq:SteinH^*H}
I_{pn}-{\mathbf H}_0^*{\mathbf H}_0=
\left(\begin{smallmatrix}0_{m\times m}&0\\
0&\Delta_{pn-m}\end{smallmatrix}\right),
\end{equation}
where $\Delta_{pn-m}$ is a $(pn-m)\times(pn-m)$ positive semi-definite
weak contraction.
\vskip 0.2cm

\mbox{Moreover, if $~p=m$ then 
$\Delta_{pn-m}$ is an orthogonal projection.}
\vskip 0.2cm

II. For $m\geq p$,
\begin{equation}\label{eq:SteinHH^*}
I_{mn}-{\mathbf H}_0{\mathbf H}_0^*
=\left(\begin{smallmatrix}0_{p\times p}&0\\
0&\Delta_{mn-p}\end{smallmatrix}\right),
\end{equation}
where $\Delta_{mn-p}$ is a $(mn-p)\times(mn-p)$ positive semi-definite
weak contraction.
\vskip 0.2cm

\mbox{Moreover, if $~p=m$ then 
$\Delta_{mn-p}$ is an orthogonal projection.}
\end{Pn}

{\bf Proof :}\quad
The structure of the right hand side of \eqref{eq:SteinH^*H} and
\eqref{eq:SteinHH^*}, is immediate from \eqref{eq:H^*H}.
All is left to verify is the spectrum of $\Delta$.
\vskip 0.2cm

Recall that the square of the Hankel singular values, are the
positive eigenvalues of
${\mathbf H}_0^*{\mathbf H}_0$ (or of ${\mathbf H}_0{\mathbf H}_0^*$)
which in turn are equal to the positive eigenvalues of the product
of the Gramians $(W_{\rm cont}W_{\rm obs})$, appeared in
\eqref{eq:gramian}, see e.g. \cite[Eq. (12.8.43)]{HSK}.
\vskip 0.2cm

From Lemma \ref{CorollaryGramians} we know that in the
rectangular case, the  spectrum of $(W_{\rm cont}W_{\rm obs})$
lies in the interval $[0,~1]$ and in the square case the non-zero
eigenvalues of $(W_{\rm cont}W_{\rm obs})$ are all ones.
\vskip 0.2cm

Finally, note that the eigenvalues of the right hand side of
\eqref{eq:SteinH^*H} and \eqref{eq:SteinHH^*} are just 1 minus
the eigenvalues of ${\mathbf H}_0^*{\mathbf H}_0$ or of
${\mathbf H}_0{\mathbf H}_0^*$, respectively.
\qed
\vskip 0.2cm

It is interesting to point out that the condition of having the
quantities in \eqref{eq:SteinH^*H} and \eqref{eq:SteinHH^*}
respectively,
\[
I_{pn}-{\mathbf H}_0^*{\mathbf H}_0\quad\quad\quad\quad\quad
I_{mn}-{\mathbf H}_0{\mathbf H}_0^*
\]
positive semidefinite (as we indeed do), commonly appears in
the context of Nehari's problem where one approximates an
anti-causal polynomial (positive powers of $z$) by a causal
rational function (no pole at infinity), see e.g.
\cite[Section 12.8]{HSK}.
\vskip 0.2cm

The following example illustrates the above results.

\begin{Ex}\label{Ex:Unitary}
{\rm
Consider part II of Example \ref{Ex:Polynomials}.\quad
As ${\mathbf H}$ in \eqref{eq:ExampHr} satisfies condition
\eqref{eq:H^*H}, from Theorem \ref{Ob:UnitaryF} it follows that
the polynomial in \eqref{eq:CoIsoEx} is 
in ${\mathcal U}$.
}
\qed
\end{Ex}
\vskip 0.2cm

The above results in Theorem \ref{Ob:UnitaryF} and in Proposition
\ref{Pn:UnitaryF} were formulated in the language of
${\mathbf H}_0^*{\mathbf H}_0$ or ${\mathbf H}_0{\mathbf H}_0^*$,
where ${\mathbf H}_0$ is a ~{\em Hankel}~ matrix. We now show
that these results can be equivalently formulated in terms
of ~{\em block-triangular Toeplitz}~ matrices. Indeed,
using ${\mathbf T}$ from \eqref{eq:T} note that
\[
\begin{matrix}
{\mathbf H}_0^*{\mathbf H}_0&=&
\left({\mathbf T}_{n, p}{\mathbf H}_0\right)^*
\left({\mathbf T}_{n, p}{\mathbf H}_0\right)&=&
\left(\begin{smallmatrix}
B_n   &  ~   &~ \\
\vdots&\ddots&~ \\
B_1   &\cdots&B_n
\end{smallmatrix}\right)^*
\left(\begin{smallmatrix}
B_n   &  ~   &~ \\
\vdots&\ddots&~ \\
B_1   &\cdots&B_n
\end{smallmatrix}\right)
\\~\\
{\mathbf H}_0{\mathbf H}_0^*&=&
\left({\mathbf H}_0{\mathbf T}_{n, m}\right)
\left({\mathbf H}_0{\mathbf T}_{n, m}\right)^*&=&
\left(\begin{smallmatrix}
B_n&\cdots&B_1   \\
  ~&\ddots&\vdots\\
~  &   ~  &B_n
\end{smallmatrix}\right)
\left(\begin{smallmatrix}
B_n&\cdots&B_1   \\
  ~&\ddots&\vdots\\
~  &   ~  &B_n
\end{smallmatrix}\right)^*
\end{matrix}
\]
Technically, this relation is well known. The problem of
finding a factorization of (block) positive semi-definite
matrix to (block) triangular Toeplitz is classical, see e.g.
\cite{Riss}. In Subsection \ref{sec:realization1} we saw
that Toeplitz operator are better established in system theory.

\section{a concluding remark}
\setcounter{equation}{0}
\label{sec:conclusion}

Although modest in size, the family of para-unitary Finite Impulse
Response systems (=(co)-isometric Laurent polynomials) is of great
interest in various fields.
\vskip 0.2cm

In Theorems \ref{Tm:RealizationU},
\ref{Tm:UnitaryFIR} and \ref{Ob:UnitaryF} we have offered three
characterizations of this family. In Corollary \ref{FIRangles} we
introduced an easy-to-use parameterization of this set.
\vskip 0.2cm

Finally, in Subsection \ref{subsec:families} we suggested six ways
(along with their combinations) to construct from a given para-unitary
FIR system, a whole family of systems, of various dimensions and
powers, all with this property.
This may raise the following open problem of going in the opposite
direction: Given a ``complicated" FIR in $\mathcal{U}$. How to
``factorize" it, following items {\bf I} to {\bf VI}, to simpler
building blocks.
\vskip 0.2cm

The nature of this work suggests
that the same set appears, under possibly different terminology,
in additional fields and further connections may be established.

\vskip 0.3cm

\begin{center}
Appendix: Construction of families of FIRs
\end{center}
\setcounter{equation}{0}
\label{section:ConstrFamilies}

For each of the items {\bf I}~ through ~{\bf VI}~
from Subsection \ref{subsec:families} we here fill-in the
following details:\\
1. Construct new polynomials from the original one.\\
2. Show that whenever the original polynomial was in
$\mathcal{U}$, so are the newly constructed polynomials.
\vskip 0.2cm

{\bf I.}\quad Reverse polynomial\\
Using ${\mathbf T}$ from \eqref{eq:T},
\[
\begin{matrix}
F_{\rm rev}(z)&:=&
z^{-1}B_n+z^{-2}B_{n-1}+~\ldots~+z^{1-n}B_2+z^{-n}B_1\\~\\
~&=&
{\mathbf Z}{\mathbf T}_{n, p}{\mathbf B}_0=
\hat{\mathbf B}{\mathbf T}_{n, m}\hat{\mathbf Z}\end{matrix}.
\]
The realtion with the corresponding Hankel matrix is
straightforward and thus omitted.
\vskip 0.2cm

{\bf II.}\quad Preserving the McMillan degree.\\
For simplicity of exposition we consider as an illustrative example
the polynomial $F(z)$ in \eqref{eq:Poly} with $n=4$ and the
parameter $q$ attaining the values $-1$ and $-2$.
\vskip 0.2cm

{\bf a.}\quad For $q=-1$ the corresponding Hankel matrix,
${\mathbf H_1}$ from
\eqref{eq:ContrR} can be, without loss of generality, extended
with a row and a column of zeros so it is $6p\times 6m$. Now
the resulting Hankel matrix may be partitioned in two forms,
\[
{\mathbf H}_1=\left(\begin{smallmatrix}
0  &B_1&B_2&B_3&B_4&0\\
B_1&B_2&B_3&B_4&0&0\\
B_2&B_3&B_4&0&0&0\\
B_3&B_4&0  &0&0&0\\
B_4&0&0&0&0&0\\
0  &0&0&0&0&0
\end{smallmatrix}\right)
\quad\quad\quad\quad
{\footnotesize\left(\begin{array}{cc|cc|cc}
0&B_1&B_2&B_3&B_4&0\\
B_1&B_2&B_3&B_4&0&0\\
\hline
B_2&B_3&B_4&0&0&0\\
B_3&B_4&0&0&0&0\\
\hline
B_4&0&0&0&0&0\\
0&0&0&0&0&0
\end{array}\right)}
\]
which produce the two polynomials in item ~{\bf II~ a}.
\vskip 0.2cm

Following Theorem \ref{Ob:UnitaryF} note that
\[
\begin{matrix}
\left(I_{6m}-{\mathbf H}_1^*{\mathbf H}_1\right)
\left(\begin{smallmatrix}
I_{2m}\\0_{4m\times 2m}\end{smallmatrix}\right)&=&0_{6m\times 2m}
&&&p\geq m\\~\\
\left(I_{2p}\quad 0_{2p\times 4p}\right)
\left(I_{6p}-{\mathbf H}_1{\mathbf H}_1^*\right)&=&0_{2p\times 6p}&&&m\geq p,
\end{matrix}
\]
so indeed the polynoial is in $\mathcal{U}$.
\vskip 0.2cm

{\bf b.}\quad Substituting in $F(z)$ in \eqref{eq:OrigPoly} $q=-2$ 
yield a Hankel matrix ${\mathbf H}_2$ with three partitionings
\[
\begin{matrix}
&
{\mathbf H}_2=
\left(\begin{smallmatrix}
0&0&B_1&B_2&B_3&B_4          \\
0&B_1&B_2&B_3&B_4&0\\
B_1&B_2&B_3&B_4&0&0\\
B_2&B_3&B_4&0&0&0\\
B_3&B_4&0&0&0&0\\
B_4&0&0&0&0&0
\end{smallmatrix}\right)&~\\~\\
{\footnotesize\left(\begin{array}{cc|cc|cc}
0&0&B_1&B_2&B_3&B_4          \\
0&B_1&B_2&B_3&B_4&0\\
\hline
B_1&B_2&B_3&B_4&0&0\\
B_2&B_3&B_4&0&0&0\\
\hline
B_3&B_4&0&0&0&0\\
B_4&0&0&0&0&0
\end{array}\right)}
&&
{\footnotesize\left(\begin{array}{ccc|ccc}
0&0&B_1&B_2&B_3&B_4          \\
0&B_1&B_2&B_3&B_4&0\\
B_1&B_2&B_3&B_4&0&0\\
\hline
B_2&B_3&B_4&0&0&0\\
B_3&B_4&0&0&0&0\\
B_4&0&0&0&0&0
\end{array}\right)}
\end{matrix}
\]
so that the three polynomials in item ~{\bf II~ b}~ are obtained.
\vskip 0.2cm

These three polynomials are in $\mathcal{U}$, as
following Theorem \ref{Ob:UnitaryF} one obtains,
\[
\begin{matrix}
\left(I_{6m}-{\mathbf H}_2^*{\mathbf H}_2\right)
\left(\begin{smallmatrix}
I_{3m}\\0_{3m\times 3m}\end{smallmatrix}\right)&=&0_{6m\times 3m}
&&&p\geq m\\~\\
\left(I_{3p}\quad 0_{3p\times 3p}\right)
\left(I_{6p}-{\mathbf H}_2{\mathbf H}_2^*\right)&=&0_{3p\times 6p}&&&m\geq p.
\end{matrix}
\]
To make the last construction more realistic take for example, $m=p=2$ and
\[
B_1={\scriptstyle\frac{1}{10}}\left(\begin{smallmatrix}3&~~3\\ 3&~~3
\end{smallmatrix}\right)\quad\quad\quad
B_2={\scriptstyle\frac{1}{5}}\left(\begin{smallmatrix}2&-2\\2&-2
\end{smallmatrix}\right)\quad\quad\quad
B_3={\scriptstyle\frac{1}{5}}\left(\begin{smallmatrix}
-2&-2\\~~~2&~~~2
\end{smallmatrix}\right),\quad\quad\quad
B_4={\scriptstyle\frac{1}{10}}\left(\begin{smallmatrix}~~~
3&-3\\-3&~~~3\end{smallmatrix}\right).
\]
\vskip 0.2cm

In items {\bf III} through {\bf VI} we produce more elaborate
structures out of a given Hankel matrix. To this end, we find it
convenient to introduce the following notation.
Let $\alpha,\beta\geq 0$ and $\eta,\delta>0$ be integers.
One can construct the following
$(\alpha+\beta+\eta)\delta\times\eta\delta$
isometry, i.e. \mbox{$U_{\rm Iso}^*U_{\rm Iso}=I_{\eta\delta}$},
\begin{equation}\label{eq:Iso}
U_{\rm Iso}(\alpha, \beta, \eta, \delta)=
I_{\eta}\otimes\left(\begin{smallmatrix}
0_{\beta\times\delta}\\I_{\delta}\\
0_{\alpha\times\delta}\end{smallmatrix}\right)=
\left(\begin{smallmatrix}
0_{\beta\times\delta} &~&~&~\\
I_{\delta}            &~&~&~\\
0_{\alpha\times\delta}&~&~&~\\
&0_{\beta\times\delta} &~&~\\
&I_{\delta }           &~&~\\
&0_{\alpha\times\delta}&~&~\\
&&\ddots&~\\
&~&~&0_{\beta\times\delta} \\
&~&~&I_{\delta}            \\
&~&~&0_{\alpha\times\delta}
\end{smallmatrix}\right)
\end{equation}
Similarly, $U_{\rm Coiso}$ is the following
$\eta\delta\times(\alpha+\beta+\eta)\delta$ co-isometry, i.e.
$U_{\rm Coiso}U_{\rm Coiso}^*=I_{\eta\delta}$,
\begin{equation}\label{eq:Coiso}
\begin{matrix}
U_{\rm Coiso}(\alpha, \beta, \eta, \delta)&=&
I_{\eta}\otimes\left(\begin{smallmatrix}0_{\delta\times\beta
}&
I_{\delta}&0_{\delta\times
\alpha
\delta}\end{smallmatrix}\right)
\\~\\~&=&
\left(\begin{smallmatrix}
0_{\delta\times\beta}&I_{\delta}&0_{\delta\times\alpha}
&~&~&~&~&~&~&~\\
&~&~&0_{\delta\times\beta}&I_{\delta}&
0_{\delta\times\alpha}&~&~&~&~\\
&~&~&~&~&~&\ddots&~&~&~\\
&~&~&~&~&~&~&0_{\delta\times\beta}&I_{\delta}&0_{\delta\times\alpha}
\end{smallmatrix}\right).
\end{matrix}
\end{equation}
Let $\rho\in\{1,~2,~\ldots~,~n\}$ be so that $\frac{n}{\rho}$ is
natural.
Substitute in $~U_{\rm Iso}$ and in $~U_{\rm Coiso}$
see \eqref{eq:Iso}, \eqref{eq:Coiso} respectively: 
$\quad\alpha=a\rho{p}$, $\beta=b\rho{p}$, with $a, b\geq 0$
integers, 
$\eta=\frac{n}{\rho}$, $\delta=\rho{p}~$
and consider the pair of products,
\begin{equation}\label{eq:UisoB}
U_{\rm Iso}{\mathbf B}=\left(\begin{smallmatrix}
0_{b\rho{p}\times m}\\B_1\\ \vdots\\~\\B_{\rho}\\
0_{(a+b)\rho{p}\times m}\\B_{\rho+1}\\\vdots\\~\\B_{2\rho}\\
0_{(a+b)\rho{p}\times m}\\
\vdots\\~\\ 0_{(a+b)\rho{p}\times m}\\
B_{n+1-\rho}\\ \vdots\\~\\B_n\\
0_{a\rho{p}\times m}
\end{smallmatrix}\right),
\end{equation}
and
\begin{equation}\label{eq:Bcoiso}
\begin{matrix}
\hat{\mathbf B}U_{\rm Coiso}&=
(0_{p\times b\rho m}~,~B_1~\cdots~B_{\rho}~,~
0_{p\times(a+b)\rho m}&~\\~\\~&
B_{\rho+1}~\cdots~B_{2\rho}~,~
0_{p\times(a+b)\rho m}~,~\cdots&
B_{n+1\rho}~\cdots~
B_n~,~0_{p\times a\rho m}).\end{matrix}
\end{equation}
Both cases yield the same $(a+b+1)np\times(a+b+1)nm$ Hankel matrix,
denoted by ${\mathbf H}(a, b, \rho)$. For example,
\begin{equation}\label{eq:Halpha1}
{\mathbf H}(0,~2,~1)=
\left(\begin{smallmatrix}
~     &~     &B_1   &~     &~     &B_2   &~     &\cdots&~     &B_n\\
~     &\Ddots&~     &~     &\Ddots&~     &\cdots&~     &\Ddots&~  \\
B_1   &~     &~     &\Ddots&~     &\cdots&~     &\Ddots&~     &~  \\
~     &~     &\Ddots&    ~ &\cdots&~     &\Ddots& ~    &~     &~  \\
~     &\Ddots&~     &\cdots&~     &\Ddots&~     &~     &~    &~  \\
B_2   &~     &\cdots&~     &\Ddots&~     &  ~   &~     &~    &~  \\
~     &~     &~     &\Ddots&~     &~     &~     &~     &~    &~  \\
~     &~     &\Ddots&~     &~     &~     &~     &~     &~    &~  \\
~     &\Ddots&~     &~     &~     &~     &~     &~     &~    &  ~\\
B_n   &~     &~     &~     &~     &~     &~     &~     &~    &~
\end{smallmatrix}\right),
\end{equation}
or
\begin{equation}\label{eq:Halpha2}
{\mathbf H}(2,~2,~2)=
\left(\begin{smallmatrix}
~      &~     &B_1   &B_2   &~     &\ldots&~     &B_{n-1}&B_n&~ \\
~      &\Ddots&\Ddots&~     &\ldots&~     &\Ddots&\Ddots&~   &~ \\
B_1    &\Ddots&~     &\ldots&~     &\Ddots&\Ddots&~     &~   &~\\
B_2    &~     &\ldots&~     &\Ddots&\Ddots&~     &~     &~   &~\\
~      &~     & ~    &\Ddots&\Ddots&~     &~     &~     &~   &~\\
~      &~     &\Ddots&\Ddots&~     &~     &~     &~     &~   &~\\
~      &\Ddots&\Ddots&~     &~     &~     &~     &~     &~   &~\\
B_{n-1}&\Ddots&~     &~     &~     &~     &~     &~     &~   &~\\
B_n    &~      &~    &~     &~     &~     &~     &~     &~   &~\\
~      &~      &~    &~     &~     &~     &~     &~     &~   &~
\end{smallmatrix}\right),
\end{equation}
or for $~n=6$
\[
{\mathbf H}(0,~1,~3)=
\left(\begin{smallmatrix}
~  & ~    &B_1   &B_2   &B_3   & ~    & ~    & ~    &B_4   &B_5   &B_6&~  \\
~  &\Ddots&\Ddots&\Ddots& ~    & ~    & ~    &\Ddots&\Ddots&\Ddots& ~ &~  \\
B_1&\Ddots&\Ddots& ~    & ~    & ~    &\Ddots&\Ddots&\Ddots&~     &~  &~  \\
B_2&\Ddots& ~    & ~    & ~    &\Ddots&\Ddots&\Ddots& ~    &~     &~  &~  \\
B_3& ~    & ~    & ~    &\Ddots&\Ddots&\Ddots& ~    & ~    & ~    &~  &~  \\
~  & ~    & ~    &\Ddots&\Ddots&\Ddots& ~    & ~    & ~    & ~    &~  &~  \\
~  & ~    &\Ddots&\Ddots&\Ddots& ~    & ~    & ~    & ~    &~     &~  &~  \\
~  &\Ddots&\Ddots&\Ddots& ~    & ~    & ~    & ~    & ~    & ~    &~  &~  \\
B_4&\Ddots&\Ddots& ~    & ~    & ~    & ~    & ~    & ~    & ~    &~  &~  \\
B_5&\Ddots& ~    & ~    & ~    & ~    & ~    & ~    & ~    & ~    &~  &~  \\
B_6& ~    & ~    & ~    & ~    & ~    & ~    & ~    & ~    & ~    &~  &~
\end{smallmatrix}\right).
\]
{\bf III.}\quad Doubling the powers.\\
First note that ${\mathbf H}(a, b, 1)$ (i.e. when $\rho=1$)
corresponds to $p\times m$-valued polynomial in \eqref{eq:z^a}
with $\gamma:=a+b+1$.
\vskip 0.2cm

As another example, the above Hankel matrix ${\mathbf H}(0,~1,~3)$ is
associated with the polynomial $F(z)$ in \eqref{eq:jumps}.
\vskip 0.2cm

{\bf IV.}\quad Rectangular polynomials.\\
Another sample of a Hankel matrix associated with
$U_{\rm Iso}{\mathbf B_o}~$ in \eqref{eq:UisoB} 
(or $~\hat{\mathbf B_o}U_{\rm Coiso}$ in \eqref{eq:Bcoiso})
is obtained when the parameters are
$~a=2$, $b=2$ $\rho=2$, i.e.
Now, multiplying ${\mathbf H}(a,~b,~\rho)$
from the ~{\em right}~ by $U_{\rm Iso}$ in
\eqref{eq:Iso} with the parameters \mbox{$\alpha=(\rho-1)m$,}
\mbox{$\beta=0$,} $\eta=\frac{n}{\rho}(a+b+1)$ and $\delta=m$
yields the following
\mbox{$(a+b+1)np\times(a+b+1)\frac{n}{\rho}m$}
Hankel matrix (here $a=2, b=2, \rho=2$)
\begin{equation}\label{eq:TallH}
{\mathbf H}(2,~2,~2)U_{\rm Iso}=
{\footnotesize\left(\begin{array}{c|c|c|c|c|c|c|c}
~      &B_1   &~     &~      &B_3   &\ldots &~      &B_{n-1}\\
~      &B_2   &~     &~      &B_4   &\ldots &~      &B_n    \\
\hline
B_1    &~     &~     &B_3    &\ldots&~      &B_{n-1}&~      \\
B_2    &~     &~     &B_4    &\ldots&~      &B_n    &~      \\
\hline
~      &~     &B_3   &\ldots &~      &B_{n-1}&~     &~      \\
~      &~     &B_4   &\ldots &~      &B_n    &~     &~      \\
\hline
~      &B_3   &\ldots&~      &B_{n-1}&~      &~     &~     \\
~      &B_4   &\ldots&~      &B_n    &~      &~     &~     \\
\hline
B_3    &\ldots&~     &B_{n-1}&~      &~      &~     &~\\
B_4    &\ldots&~     &B_n    &~      &~      &~     &~\\
\hline
\cdots &\vdots&\vdots&~      &~      &~      &~     &~\\
\hline
~      &B_{n-1}&~    &~      &~      &~      &~     &~\\
~      &B_n    &~    &~      &~      &~      &~     &~\\
\hline
B_{n-1}&~      &~    &~      &~      &~      &~     &~\\
B_n    &~      &     &~      &~      &~      &~     &~\\
\hline
~      &~      &~    &~      &~      &~      &~     &~\\
  ~    &~      &~    &~      &~      &~      &~     &~
\end{array}\right)},
\end{equation}
which corresponds to the $p\rho\times m$-valued
polynomial in {\bf IV~a.}
\vskip 0.2cm

Similarly, multiplying ${\mathbf H}(a,~b,~\rho)$
from the ~{\em left}~ by $U_{\rm Coiso}$ in
\eqref{eq:Coiso} with the parameters $~\alpha=(\rho-1)p$,
$\beta=0$, $\eta=\frac{n}{\rho}(a+b+1)$ and $\delta=p~$
yields the following
\mbox{$(a+b+1)\frac{n}{\rho}p\times(a+b+1)nm$}
Hankel matrix (here $a=2, b=2, \rho=2$)
\begin{equation}\label{eq:FatH}
U_{\rm Coiso}{\mathbf H}(2,~2,~2)=
{\footnotesize\left(\begin{array}{c|c|c|c|c|c}
~  &B_1\quad B_2&~  &\ldots&B_{n-1}\quad B_n&      ~\\
\hline
B_1\quad B_2&~   &~  &\ldots&~      &~\\
\hline
~      &~  &B_3\quad B_4&\ldots&~      &~\\
\hline
~  &B_3\quad B_4 &~  &\ldots &~      &~\\
\hline
B_3\quad B_4&~   &~  &\ldots  &~      &~\\
\hline
\vdots&\vdots&\vdots&
\cdots&  ~   &~\\
\hline
~      &B_{n-1}\quad B_n&~&  ~    &~&~  \\
\hline
B_{n-1}\quad B_n    &~&~&  ~   &~&~\\
\hline
~&~&~&~     &~&~
\end{array}\right)},
\end{equation}
which corresponds to the $p\times\rho{m}$-valued
polynomial in {\bf IV~b.}
\vskip 0.2cm

It is easy to verify that if ${\mathbf H_o}$ satisfies
the first line in \eqref{eq:H^*H} then so do all Hankel
matrices of the form ${\mathbf H}_0U_{\rm Iso}$
\eqref{eq:TallH}.
\vskip 0.2cm

Similarly, if ${\mathbf H}_0$ satisfies the
second line in \eqref{eq:H^*H}
then so do all Hankel matrices of the form 
$U_{\rm Coiso}{\mathbf H}_0$ \eqref{eq:FatH}.
\vskip 0.2cm

In the sequel, we shall adjust our previous notation
in \eqref{eq:ContrR} of the Hankel matrix associated with
\[
F(z)=z^{-(1+\eta)}B_1+~\ldots~+z^{-(n+\eta)}B_n
\quad\quad\quad\quad\eta\geq 0,
\]
to $~
{\mathbf H}_{\mathbf{B}, n, \eta}
$
(in \eqref{eq:ContrR} the subscricts $B$ and $n$ were omitted, as so far 
they were evident from the context). For example, with the polynomial
\[
z^{-(1+\eta)}C_1+~\ldots~+z^{-(l+\eta)}C_l
\quad\quad\quad\quad\eta=0,~1,~2\ldots
\]
one can associate the $(l+\eta)p_c\times(l+\eta)m_c$ Hankel matrix,
\[
{\mathbf H}_{\mathbf{C}, l, \eta}=\left(\begin{smallmatrix}
~     &~     &C_1   &\cdots&C_l\\
~     &\Ddots&   ~  &\Ddots&   \\
C_1   &\cdots&C_l   & ~    & \\
\vdots&\Ddots&~     & ~    & \\
C_l&~ & ~    & ~    & ~    & 
\end{smallmatrix}\right).
\]
{\bf V.}\quad Composition of polynomials.\\
With the pair of poynomilas in \eqref{eq:PairPoly} one can
associated the Hankel matrices ${\mathbf H}_{\mathbf{B}, n, 0}$
and ${\mathbf H}_{\mathbf{C}, l, 0}$, which are of dimensions
$np_b\times nm_b$ and $lp_c\times nm_c$,
respectively.
\vskip 0.2cm

Out of this pair, one can construct (at least) the ~{\em three}~
following Hankel matrices, all of the form
${\mathbf H}_{\mathbf{D}, n, 0}~$:\\
{\bf a.}\quad A $~n(p_b+p_c)\times n(m_b+m_c)$ Hankel matrix
\[
{\footnotesize\left(\begin{array}{cc|cc|c|cc|cc|c|cc}
B_1&~  &B_2&~  &\ldots&B_l    &~  &B_{l+1}&~&\ldots&B_n&~\\
~  &C_1&~  &C_2&\ldots&~      &C_l&~      &0_{p_c\times m_c}
&\ldots&~  &0_{p_c\times m_c}\\
\hline
B_2&~  &B_3&~  &\ldots&B_{l+1}&~  &~      &~&~&~  &~\\
~  &C_2&~  &C_3&\ldots&~      &0_{p_c\times m_c}&~&~&~&~  &~\\
\hline
B_3&~  &B_4&~  &\ldots&~      &~  &~      &~&~&~  &~\\
~  &C_3&~  &C_4&\ldots&~      &~  &~      &~&~&~  &~\\
\hline
\vdots&\vdots&\vdots&\vdots&~&~&
~&~&~&~&   ~  &   ~  \\
\hline
B_l&~  &B_{l+1}&~&\ldots&~      &~  &~      &~&  ~   &~  &~\\
~  &C_l&~      &0_{p_c\times m_c}&\ldots&~      &~  &~      &~&   ~  &~  &~\\
\hline
B_{l+1}&~&~&~&  ~   &~      &~  &~      &~&  ~   &~  &~\\
~      &0_{p_c\times m_c}&~&~&  ~   &~      &~  &~      &~&   ~  &~  &~\\
\hline
\vdots&\vdots&~&   ~  &  ~   &   ~  &  ~   &   ~  &~&   ~  &  ~   &  ~   \\
\hline
B_n&~&~&~&  ~&~      &~  &~      &~&   ~  &~  &~\\
~&0_{p_c\times m_c}&~&~&   ~  &~      &~  &~      &~&   ~  &~  &~
\end{array}\right)}
\]
or another
$~n(p_b+p_c)\times n(m_b+m_c)$ Hankel matrix
\[
{\footnotesize\left(\begin{array}{cc|cc|c|cc|cc|c|cc}
~  &B_1&~  &B_2&\ldots&~                &B_l     &~                &B_{l+1}&\ldots&~&B_n\\
C_1&~  &C_2&~  &\ldots&C_l              &~       &0_{p_c\times m_c}&~      &\ldots&0_{p_c\times m_c}&0\\
\hline
~  &B_2&~  &B_3&\ldots&~                &B_{l+1} &~                &~      &~     &~&~  \\
C_2&~  &C_3&~  &\ldots&0_{p_c\times m_c}&~       &~                &~      &~     &~&~ \\
\hline
~  &B_3&~  &B_4&\ldots&~                &~       &~                &~      &~     &~&~\\
C_3&~  &C_4&~  &\ldots&~                &~       &~                &~      &~     &~&~\\
\hline
\vdots&\vdots&\vdots&\vdots&~&~&
~&~&~&~&   ~  &   ~  \\
\hline
~  &B_l&~   &B_{l+1}&~&~      &~  &~      &~&  ~   &~  &~\\
C_l&~  &0_{p_c\times m_c}&~&~      &~  &~      &~&   ~  &~  &~\\
\hline
~&B_{l+1}&~&~&  ~   &~      &~  &~      &~&  ~   &~  &~\\
0_{p_c\times m_c}&~&~&  ~   &~      &~  &~      &~&   ~  &~  &~\\
\hline
\vdots&\vdots&~&   ~  &  ~   &   ~  &  ~   &   ~  &~&   ~  &  ~   &  ~   \\
\hline
~&B_n&~&~&  ~&~      &~  &~      &~&   ~  &~  &~\\
0_{p_c\times m_c}&~&~&   ~  &~      &~  &~      &~&   ~  &~  &~
\end{array}\right)}
\]
\vskip 0.2cm

{\bf b.}\quad For $m_c\geq m_b~$ a 
$~n(p_b+p_c)\times nm_c$ Hankel matrix
\[
{\footnotesize\left(\begin{array}{c|c|c|c|c|c}
\sqrt{\alpha}B_1~~0_{p_b\times(m_c-m_b)}
&~&\sqrt{\alpha}B_l~~0_{p_b\times(m_c-m_b)}
&\sqrt{\alpha}B_{l+1}~~0_{p_b\times(m_c-m_b)}
&~
&\sqrt{\alpha}B_n~~0_{p_b\times(m_c-m_b)}\\
\sqrt{1-\alpha}C_1&~&\sqrt{1-\alpha}C_l&0_{p_c\times m_c}&~
&0_{p_c\times m_c}\\
\hline
\vdots&\Ddots&~&~&\Ddots&~\\
\hline
\sqrt{\alpha}B_l~~0_{p_b\times(m_c-m_b)}&~&~&
\sqrt{\alpha}B_n~~0_{p_b\times(m_c-m_b)}&~&~\\
\sqrt{1-\alpha}C_l&~&~&0_{p_c\times m_c}&~&~\\
\hline
\sqrt{\alpha}B_{l+1}~~0_{p_b\times(m_c-m_b)}&~&
\sqrt{\alpha}B_n~~0_{p_b\times(m_c-m_b)}
&~                   &~\\
0_{p_c\times m_c}     &~                 &0_{p_c\times m_c}                   &~&~\\
\hline
\vdots&\Ddots& ~    &  ~   &   ~  \\
\hline
\sqrt{\alpha}B_n~~0_{p_b\times(m_c-m_b)}&~ &~                   &~&  ~  \\
0_{p_c\times m_c}     &~                 &  ~   &~                   &~&~
\end{array}\right)}
\]
\vskip 0.2cm

{\bf c.}\quad For $p_b\geq p_c$ a $np_b\times(m_b+m_c)n$ Hankel matrix
\[
{\footnotesize\left(\begin{array}{cc|c|cc|cc|c|cc}
\begin{smallmatrix}\sqrt{\alpha}B_1\\ 0_{(p_b-p_c)\times m_b}\end{smallmatrix}
&\sqrt{1-\alpha}C_1&~&
\begin{smallmatrix}\sqrt{\alpha}B_l\\ 0_{(p_b-p_c)\times m_b}\end{smallmatrix}
&\sqrt{1-\alpha}C_l&
\begin{smallmatrix}\sqrt{\alpha}B_{l+1}\\ 0_{(p_b-p_c)\times m_b}\end{smallmatrix}&
0_{p_b\times m_c}&~&
\begin{smallmatrix}\sqrt{\alpha}B_n\\ 0_{(p_b-p_c)\times m_b}\end{smallmatrix}
&0_{p_b\times m_c}\\
\hline
\vdots&~&\Ddots&~&~&~&~&\Ddots&~&~\\
\hline
\begin{smallmatrix}\sqrt{\alpha}B_l\\ 0_{(p_b-p_c)\times m_b}\end{smallmatrix}
&\sqrt{1-\alpha}C_l&~&~&~&
\begin{smallmatrix}\sqrt{\alpha}B_n\\ 0_{(p_b-p_c)\times m_b}\end{smallmatrix}
&0_{p_b\times m_c}&~&~\\
\hline
\begin{smallmatrix}\sqrt{\alpha}B_{l+1}\\ 0_{(p_b-p_c)\times m_b}\end{smallmatrix}
&0_{p_b\times m_c}&~&
\begin{smallmatrix}\sqrt{\alpha}B_n\\ 0_{(p_b-p_c)\times m_b}\end{smallmatrix}
&0_{p_b\times m_c}&~&~&~&~&~\\
\hline
\vdots&\vdots&\Ddots&
&&~&  ~  &   ~  &~\\
\hline
\begin{smallmatrix}\sqrt{\alpha}B_n\\ 0_{(p_b-p_c)\times m_b}\end{smallmatrix}
&0_{p_b\times m_c}&~&~&~&~&~&~&~&~
\end{array}\right)}
\]

{\bf VI.}\quad Product of polynomials\\
Recall that out of
\[
\begin{matrix}
F_b(z)&=&z^{-1}B_1+\ldots~z^{-n}B_n&&p_b\times\rho
\\~\\
F_c(z)&=&z^{-1}C_1+\ldots~z^{-l}C_l&&\rho
\times m_c
\end{matrix}
\]
the following $~p_b\times m_c$-valued polynomial was obtained 
\begin{equation}\label{eq:Fd}
F_d(z):=F_b(z)F_c(z)=z^{-1}\left(z^{-1}D_1+~\ldots~+z^{-(n+l-1)}D_{n+l-1}\right)~.
\end{equation}
where the coefficients $D_1~,~\ldots~,~D_{n+l-1}~$ were
explicitely given in \eqref{eq:CoeffD}.
\vskip 0.2cm

Expressing, \eqref{eq:CoeffD} in terms of corresponding
Hankel matrices yields
\begin{equation}\label{eq:Hd}
{\mathbf H}_{\mathbf{D}, n+l-1, 1}=
{\mathbf H}_{\mathbf{B}, n, l}{\mathbf T}_{m+l, \rho}
{\mathbf H}_{\mathbf{C}, l, n}~,
\end{equation}
where the Hankel matrices
${\mathbf H}_{\mathbf{B}, n, l}$,
${\mathbf H}_{\mathbf{C}, l, n}~$ and
$~{\mathbf H}_{\mathbf{D}, n+l-1, 1}~$
are $(n+l)p_b\times(n+l)\rho$, $(n+l)\rho\times(n+l)m_c~$ and
$~(n+l)p_b\times(n+l)m_c~$ respectively, while
$~{\mathbf T}_{m+l, \rho}$ is the permutation matrix as in \eqref{eq:T}.
\vskip 0.2cm

Next, to establish the fact that $F_d(z)$ is in $~\mathcal{U}$,
we go through the following steps.
\vskip 0.2cm

First, from \eqref{eq:Hd} note that
\[
{\mathbf H}_{\mathbf{D}, n+l-1, 1}^*{\mathbf H}_{\mathbf{D}, n+l-1 , 1}
={\mathbf H}_{\mathbf{C}, l, n}^*{\mathbf T}_{m+l, \rho}
{\mathbf H}_{\mathbf{B}, n, l}^*{\mathbf H}_{\mathbf{B}, n, l}
{\mathbf T}_{m+l, \rho}{\mathbf H}_{\mathbf{C}, l, n}~.
\]
Assuming now that $~F_b\in\mathcal{U}$, it follows from
\eqref{eq:H^*H} that
\[
\begin{matrix}
{\mathbf H}_{\mathbf{D}, n+l-1, 1}^*
{\mathbf H}_{\mathbf{D}, n+l-1, 1}&=
&{\mathbf H}_{\mathbf{C}, l , n}^*{\mathbf T}_{m+l, \rho}\cdot
{\rm diag}\{I_{(l+1)\rho}\quad\Delta_{(n-1)\rho}\}\cdot
{\mathbf T}_{m+l, \rho}{\mathbf H}_{\mathbf{C}, l,n}\\~\\~&=&
{\mathbf H}_{\mathbf{C}, l, n}^*\cdot
{\rm diag}\{\Delta_{(n-1)\rho}\quad I_{(l+1)\rho}\}\cdot
{\mathbf H}_{\mathbf{C}, l, n}~,
\end{matrix}
\]
where $\Delta_{(n-1)\rho}$ is $~(n-1)\rho\times(n-1)\rho~$
positive semi-definite (weak) contraction.
\vskip 0.2cm

Assuming now that also $~F_c\in\mathcal{U}$, (carefully
following the dimensions) it follows
from \eqref{eq:H^*H} that
\[
{\mathbf H}_{\mathbf{D}, n+l-1, 1}^*
{\mathbf H}_{\mathbf{D}, n+l-1, 1}=
{\rm diag}\{I_{2p_b}\quad\hat{\Delta}_{(n+l-2)p_b}\},
\]
where $\hat{\Delta}_{(n+l-2)p_b}$ is a $~(n+l-1)p_b\times(n+l-1)p_b~$ 
positive semi-definite (weak) contraction, this part is estabilished
and indeed $F_d\in\mathcal{U}$.
\vskip 0.2cm

Showing the relation for
${\mathbf H}_{\mathbf{D}, n+l-1, 1}
{\mathbf H}_{\mathbf{D}, n+l-1, 1}^*$
is quite similar and thus omitted.
\qed
\end{document}